\documentclass[11pt]{article}

\usepackage{amsmath,amssymb}
\usepackage{theorem}

\makeatletter
\@addtoreset{equation}{section}
\makeatother

\newtheorem{theorem}{Theorem}[section]
\newtheorem{proposition}[theorem]{Proposition}
\newtheorem{lemma}[theorem]{Lemma}

\newtheorem{definition}[theorem]{Definition}
{\theorembodyfont{\rmfamily}

\newtheorem{remark}[theorem]{Remark}}

\title{\Huge\bf \boldmath $U_q\left(\mathfrak{sl}_2\right)$-symmetries of the quantum disc: a complete list}
\author{\Large Sergey D. Sinel'shchikov
\\ \it Mathematics Division
\\ \it B. Verkin Institute for Low Temperature Physics and Engineering
\\ \it of the National Academy of Sciences of Ukraine
\\ 47 Nauky Ave., 61103 Kharkiv, Ukraine
\\ E-mail: sinelshchikov@ilt.kharkov.ua
}
\date{}

\begin{document}
\large

\maketitle

\begin{abstract}
This work presents a classification of $U_q(\mathfrak{sl}_2)$-symmetries on the quantum disc. The principal invariant of such classification, the grading jump, is introduced. It turns out that, under the present subjects, the grading jump can take only 3 values: $0$, $1$, $-1$. The subcollection of the complete collection of symmetries is extracted in such a way that the selected symmetries satisfy certain compatibility condition for involutions.
\end{abstract}

\begin{quotation}\small
\vspace{3ex} {\it Key words}: quantum universal enveloping algebra; Hopf
algebra; quantum disc; quantum symmetry, grading jump, weight, involution

\end{quotation}

\section{Introduction}

An essential idea in studying quantum algebras is to consider them together with a certain collection of `quantum symmetries'. Normally such pairs of subjects were treated as $q$-analogs for actions of Lie groups on their homogeneous spaces.

Initially a single distinguished symmetry on the quantum plane has been considered (the original term was `the structure of $U_q(\mathfrak{sl}_2)$-module algebra on the quantum plane', see, e.g., \cite{K}); one had also a similar distinguished such structure on the quantum disc \cite{V}, just one more simplest quantum algebra to be considered in this work.

A complete list of $U_q(\mathfrak{sl}_2)$-symmetries on the quantum plane has been described in \cite{DS}. This initial result has been extended to certain quantum spaces of higher dimension, along with the related actions of $U_q(\mathfrak{sl}_n)$ by symmetries \cite{DHL}.

Another reasonable extension of the results of \cite{DS} is presented in \cite{S1,S2}, where the standard (polynomial algebra of) quantum plane is embedded into a larger quantum algebra of Laurent polynomials on the quantum plane. The latter algebra, while retaining all the symmetries of the standard quantum plane, appears to be much more symmetric, with rather extended classification list of symmetries.

The purpose of this paper is to produce a complete list of $U_q(\mathfrak{sl}_2)$-symmetries on the quantum disc $\operatorname{Pol}(\mathbb{D})_q$. Our initial assumption is that the algebra $\operatorname{Pol}(\mathbb{D})_q$ carries no involution. This was made implicit within the principal part of the research, just to obtain the utmost list of the symmetries. This list is given here in Table 1 for the reader's convenience; the notation involved therein can be found in the rest of the text.

\textbf{Table 1.}
$$
\begin{array}{||c|c|c||}
\hline\hline
\begin{array}{c}\textbf{Series}\\ \textbf{names}
\end{array} &
\begin{array}{c}\textbf{Weight}\\ \textbf{constants}
\end{array} &
\textbf{Action of $\mathsf{e}$ and $\mathsf{f}$}
\\ \hline\hline
\mathbf{(0+)} &
\begin{aligned}
\mathsf{k}(z) &=z &
\\ \mathsf{k}(z^*) &=z^* &
\end{aligned} &
\begin{aligned}
\mathsf{e}(z) &=\mathsf{e}(z^*)=0 &
\\ \mathsf{f}(z) &=\mathsf{f}(z^*)=0 &
\end{aligned}
\\ \hline
\mathbf{(0-)} &
\begin{aligned}
\mathsf{k}(z) &=-z &
\\ \mathsf{k}(z^*) &=-z^* &
\end{aligned} &
\begin{aligned}
\mathsf{e}(z) &=\mathsf{e}(z^*)=0 &
\\ \mathsf{f}(z) &=\mathsf{f}(z^*)=0 &
\end{aligned}
\\ \hline
\mathbf{(1a)} &
\begin{aligned}
\mathsf{k}(z) &=q^2z
\\ \mathsf{k}(z^*) &=q^{-2}z^*
\end{aligned} &
\begin{aligned}
\mathsf{e}(y) &=q^{-1}b_0^{-1}zy &\mathsf{f}(y) &=\left(b_0y+b_1y^2\right)z^*&
\\ \mathsf{e}(z) &=qb_0^{-1}z^2 &\mathsf{f}(z) &=-b_0-b_1y^2&
\\ \mathsf{e}(z^*) &=-q^{-1}b_0^{-1} &\mathsf{f}(z^*) &=q^2b_0z^{*2}&
\\ b_0,b_1 &\in\mathbb{C},\qquad b_0\ne 0&&&
\end{aligned}
\\ \hline
\mathbf{(1b)} &
\begin{aligned}
\mathsf{k}(z) &=q^2z
\\ \mathsf{k}(z^*) &=q^{-2}z^*
\end{aligned} &
\begin{aligned}
\mathsf{e}(y) &=z\left(a_0y+a_1y^2\right) &\mathsf{f}(y) &=q^{-1}a_0^{-1}yz^*&
\\ \mathsf{e}(z) &=q^2a_0z^2 &\mathsf{f}(z) &=-q^{-1}a_0^{-1}&
\\ \mathsf{e}(z^*) &=-a_0-a_1y^2 &\mathsf{f}(z^*) &=qa_0^{-1}z^{*2}&
\\ a_0,a_1 &\in\mathbb{C},\qquad a_0\ne 0&&&
\end{aligned}
\\ \hline
\mathbf{(-1a)} &
\begin{aligned}
\mathsf{k}(z) &=q^{-2}z
\\ \mathsf{k}(z^*) &=q^2z^*
\end{aligned} &
\begin{aligned}
\mathsf{e}(y) &=-qb_1^{-1}z^* &\mathsf{f}(y) &=z(b_0+b_1y)&
\\ \mathsf{e}(z) &=q^{-1}b_1^{-1} &\mathsf{f}(z) &=-q^2b_1z^2&
\\ \mathsf{e}(z^*) &=0 &\!\!\mathsf{f}(z^*) &=-q^{-2}b_0+b_1-\left(1+q^{-2}\right)b_1y\!\! &
\\ b_0,b_1 &\in\mathbb{C},\qquad b_1\ne 0&&&
\end{aligned}
\\ \hline
\mathbf{(-1b)} &
\begin{aligned}
\mathsf{k}(z) &=q^{-2}z
\\ \mathsf{k}(z^*) &=q^2z^*
\end{aligned} &
\begin{aligned}
\mathsf{e}(y) &=(a_0+a_1y)z^* &\mathsf{f}(y) &=-qa_1^{-1}z&
\\ \mathsf{e}(z) &=-q^{-2}a_0+a_1-\left(1+q^{-2}\right)a_1y &\mathsf{f}(z) &=0&
\\ \mathsf{e}(z^*) &=-q^2a_1z^{*2} &\mathsf{f}(z^*) &=q^{-1}a_1^{-1}&
\\ a_0,a_1 &\in\mathbb{C},\qquad a_1\ne 0&&&
\end{aligned}
\\ \hline\hline
\end{array}
$$

After that, in the last Section \ref{invs}, the subcollections of symmetries are extracted, which, under various additional assumptions on $q$ and choices of involution on $U_q(\mathfrak{sl}_2)$, are subject to a speci1al compatibility assumption on involutions.

The outline of this paper is as follows. Section \ref{prel} contains some preliminary material: definitions, notations, some well known and obvious facts. Section \ref{trl} describes the trivial series $\mathbf{(0+)}$ and $\mathbf{(0-)}$ of $U_q(\mathfrak{sl}_2)$-symmetries on $\operatorname{Pol}(\mathbb{D})_q$, together with the principal invariant of the symmetries in question, the grading jump $\operatorname{GJ}$. Section \ref{GJ+} presents a description of symmetries with $\operatorname{GJ}>0$ and demonstrates that in fact the only possible value of $\operatorname{GJ}$ for such symmetries is $\operatorname{GJ}=1$. Section \ref{GJ-} investigates the case $\operatorname{GJ}<0$; similarly, it turns out that such symmetries exist only in the case $\operatorname{GJ}=-1$. Finally, Section \ref{invs} extracts the subcollections of those $U_q(\mathfrak{sl}_2)$-symmetries on $\operatorname{Pol}(\mathbb{D})_q$ which respect involutions in the above algebras.

\section{Preliminaries}\label{prel}

We start with recalling the general definitions. Let $H$ be a Hopf algebra whose comultiplication is $\Delta$, counit is $\varepsilon$, and antipode is $S$ \cite{abe}. Consider also a unital algebra $A$ whose unit is $\mathbf{1}$. The Sweedler sigma-notation related to the comultiplication $\Delta(h)=\sum\limits_{(h)}h_{(1)}\otimes h_{(2)}$ as in \cite{sweedler} is used below. In what follows, $\mathbb{C}$ is assumed to be the ground field.

\begin{definition}\label{symdef}
By a structure of $H$-module algebra on $A$ (to be referred to as an $H$-symmetry for the sake of brevity, or even merely a symmetry if $H$ and $A$ are completely determined by the context) we mean a homomorphism of algebras $\pi\colon H\to\operatorname{End}_\mathbb{C}A$ such that
\begin{description}
\item[(i)] $\pi(h)(ab)=\sum\limits_{(h)}\pi\left(h_{(1)}\right)(a)\cdot
    \pi\left(h_{(2)}\right)(b)$ for all $h\in H$, $a,b\in A$;
\item[(ii)] $\pi(h)(\mathbf{1})=\varepsilon(h)\mathbf{1}$ for all $h\in H$.
\end{description}
The symmetries $\pi_1$, $\pi_2$ are said to be isomorphic if there exists an automorphism $\Psi$ of the algebra $A$ such that $\Psi\pi_1(h)\Psi^{-1}=\pi_2(h)$ for all $h\in H$.
\end{definition}

Throughout the paper we assume that $q\in\mathbb{C}{\setminus}\{0\}$ is not a root of $1$ ($q^n\ne 1$ for all non-zero integers $n$).

The quantum disc \cite{KL,NN,SSV1} is a unital algebra $\operatorname{Pol}(\mathbb{D})_q$ generated by $z$, $z^*$ subject to the relation
\begin{equation}\label{zz*}
zz^*=q^2z^*z+1-q^2.
\end{equation}
Certainly this is a $*$-algebra under the natural involution $z\mapsto z^*$. However, the principal purpose of this paper is to produce a complete list of $U_q\left(\mathfrak{sl}_2\right)$-symmetries on $\operatorname{Pol}(\mathbb{D})_q$, with the latter being considered as an algebra without involution, so that $z^*$ is treated as a single symbol. This is our approach before the last Section \ref{invs}, in which the details related to involutions are expounded.

We use the obvious grading on $\operatorname{Pol}(\mathbb{D})_q$ given by
$$
\operatorname{Pol}(\mathbb{D})_q=\mathop{\oplus}\limits_{k=-\infty}^\infty
\mathcal{A}_k,\qquad\text{\ with\ }\qquad\mathcal{A}_k=\text{linear span of}\left\{\left.z^iz^{*j}\right|\:i-j=k\right\}.
$$
The algebra $\operatorname{Pol}(\mathbb{D})_q$ contains an element $y=1-zz^*\in\mathcal{A}_0$, which satisfies the following quasicommutation relations
\begin{align}
yz &=q^{-2}zy,\label{yz}
\\ yz^* &=q^2z^*y.\label{yz*}
\end{align}
The general form of an element of $\mathcal{A}_k$ is $z^k\varphi(y)$ for $k\ge 0$, and $\psi(y)(z^*)^{-k}$ for $k<0$, which is an easy consequence of \eqref{zz*}. It is also worth mentioning a closely related and quite obvious fact that $\operatorname{Pol}(\mathbb{D})_q$ is a domain (no zero divisors).

The quantum universal enveloping algebra $U_q(\mathfrak{sl}_2)$ \cite{K,KS} is a unital associative algebra defined by its (Chevalley) generators $\mathsf{k}$, $\mathsf{k}^{-1}$, $\mathsf{e}$, $\mathsf{f}$, and the relations
\begin{gather}
\mathsf{k}^{-1}\mathsf{k} =\mathbf{1},\qquad\mathsf{kk}^{-1}=\mathbf{1},\notag\\
\mathsf{ke} =q^2\mathsf{ek},\label{ke}\\
\mathsf{kf} =q^{-2}\mathsf{fk},\label{kf}\\
\mathsf{ef}-\mathsf{fe} =\frac{\mathsf{k}-\mathsf{k}^{-1}}{q-q^{-1}}. \label{effe}
\end{gather}

The standard Hopf algebra structure on $U_q(\mathfrak{sl}_2)$ is determined by the comultiplication $\Delta$, the counit $\boldsymbol\varepsilon$, and the antipode $\mathsf{S}$ as follows
\begin{align}
\Delta(\mathsf{k}) &=\mathsf{k}\otimes\mathsf{k}, && && \label{k0}
\\ \Delta(\mathsf{e}) &=\mathbf{1}\otimes\mathsf{e}+ \mathsf{e}\otimes\mathsf{k}, && && \label{def}
\\ \Delta(\mathsf{f}) &=\mathsf{f}\otimes\mathbf{1}+\mathsf{k}^{-1}\otimes\mathsf{f}, && &&  \label{def1}
\\ \mathsf{S}(\mathsf{k}) &=\mathsf{k}^{-1}, & \mathsf{S}(\mathsf{e}) &=-\mathsf{ek}^{-1}, & \mathsf{S}(\mathsf{f}) &=-\mathsf{kf}, \notag
\\ \boldsymbol{\varepsilon}(\mathsf{k}) &=\mathbf{1}, & \boldsymbol{\varepsilon}(\mathsf{e}) &=\boldsymbol{\varepsilon}(\mathsf{f})=0.\qquad && \notag
\end{align}

Here and in what follows we describe the (series of) $U_q(\mathfrak{sl}_2)$-symmetries on $\operatorname{Pol}(\mathbb{D})_q$ via determining an action of the distinguished generators of $U_q(\mathfrak{sl}_2)$ on the generators of $\operatorname{Pol}(\mathbb{D})_q$. To derive the associated $U_q(\mathfrak{sl}_2)$-symmetry, we first extend the action to monomials (both in~$U_q(\mathfrak{sl}_2)$ and in $\operatorname{Pol}(\mathbb{D})_q$) using
\begin{gather*}
(\mathsf{ab})u \overset{\rm def}{=}\mathsf{a}(\mathsf{b}u), \qquad \mathsf{a},\mathsf{b} \in U_q(\mathfrak{sl}_2), \qquad u \in\operatorname{Pol}(\mathbb{D})_q,
\\ \mathsf{a}(uv) \overset{\rm def}{=} \sum_{(\mathsf{a})}\left(\mathsf{a}_{(1)}u\right)\cdot
\left(\mathsf{a}_{(2)}v\right),\qquad\mathsf{a}\in U_q(\mathfrak{sl}_2), \qquad u,v \in\operatorname{Pol}(\mathbb{D})_q,
\end{gather*}
and then extend by linearity to the entire algebras $U_q(\mathfrak{sl}_2)$ and $\operatorname{Pol}(\mathbb{D})_q$, using
\begin{gather*}
\mathsf{a}(u+v)=\mathsf{a}u+\mathsf{a}v,\qquad (\mathsf{a}+\mathsf{b})u=\mathsf{a}u+\mathsf{b}u,\\
\mathbf{1}u=u,\qquad \mathsf{a}\mathbf{1}=\boldsymbol{\varepsilon}(\mathsf{a})\mathbf{1}, \qquad\mathsf{a},\mathsf{b}\in U_q(\mathfrak{sl}_2),\qquad u,v\in\operatorname{Pol}(\mathbb{D})_q.
\end{gather*}
Such extension determines a well-defined action of $U_q(\mathfrak{sl}_2)$ on $\operatorname{Pol}(\mathbb{D})_q$ if and only if everything passes through the relations in $U_q(\mathfrak{sl}_2)$ and in $\operatorname{Pol}(\mathbb{D})_q$. To verify this, one has to apply every generator of $U_q(\mathfrak{sl}_2)$ to each relation in $\operatorname{Pol}(\mathbb{D})_q$, and then every relation in $U_q(\mathfrak{sl}_2)$ to each generator of $\operatorname{Pol}(\mathbb{D})_q$. This is to be done in each specific case, and normally such verification is left to the reader.

Given a $U_q(\mathfrak{sl}_2)$-symmetry on $\operatorname{Pol}(\mathbb{D})_q$, the generator $\mathsf{k}$ acts via an automorphism of $\operatorname{Pol}(\mathbb{D})_q$, as one can readily deduce from invertibility of $\mathsf{k}$, Definition \ref{symdef}(i) and \eqref{k0}.

A description of automorphisms of the algebra $\operatorname{Pol}(\mathbb{D})_q$ is due to J. Alev and M. Chamarie.

\begin{proposition}[\cite{ale/cha}, Proposition 1.4.4(i)]\label{auts}
Let $\Psi$ be an automorphism of $\operatorname{Pol}(\mathbb{D})_q$, then there exists a non-zero constant $\alpha$ such that
$$\Psi: z\mapsto\alpha z,\qquad z^*\mapsto\alpha^{-1}z^*.$$
\end{proposition}

This automorphism is well-defined on the entire algebra $\operatorname{Pol}(\mathbb{D})_q$, because the ideal of relations generated by \eqref{zz*} is $\Psi$-invariant.

It follows from Proposition \ref{auts} that, given a symmetry $\pi$, the action of $\mathsf{k}$ is determined completely on the generators of $\operatorname{Pol}(\mathbb{D})_q$ as follows
\begin{equation}\label{wc}
\pi(\mathsf{k})(z)=\alpha z,\qquad\pi(\mathsf{k})(z^*)=\alpha^{-1}z^*
\end{equation}
for some {\it weight constant} $\alpha\in\mathbb{C}\setminus\{0\}$. Therefore every monomial $z^iz^{*j}\in\operatorname{Pol}(\mathbb{D})_q$ is an eigenvector for $\pi(\mathsf{k})$ (a {\it weight vector}), and the associated eigenvalue $\alpha^{i-j}$ will be referred to as the {\it weight} of this monomial, to be written as $\mathbf{wt}\left(z^iz^{*j}\right)=\alpha^{i-j}$.

\begin{remark}
Observe that $\mathbf{wt}(y)=1$ and, more generally, for $u\in\mathcal{A}_0$ one has $\mathbf{wt}(u)=1$. This already implies that $\mathbf{wt}$ is constant on every homogeneous component $\mathcal{A}_k$. It is convenient to consider, instead of monomials of $z$, $z^*$, the weight vectors in the general form $z^k\varphi(y)$ and $\psi(y)(z^*)^k$, with $k\ge 0$ and $\varphi$, $\psi$ polynomials. Here, $\mathbf{wt}\left(z^k\varphi(y)\right)=k$ and $\mathbf{wt}\left(\psi(y)(z^*)^k\right)=-k$.
\end{remark}

\section{The trivial series of symmetries. The grading jump (\boldmath$\operatorname{GJ}$) related to a symmetry}\label{trl}

We start with the simplest case in which the operators $\pi(\mathsf{e})$ and $\pi(\mathsf{f})$ are identically zero.

\begin{lemma}\label{iz}
Let $\pi$ be a $U_q(\mathfrak{sl}_2)$-symmetry on $\operatorname{Pol}(\mathbb{D})_q$. The following properties of $\pi$ are equivalent:
\begin{description}
\item[(i)] the weight constant $\alpha\in\{-1;1\}$;

\item[(ii)] $\pi(\mathsf{e})$ is the identically zero operator on $\operatorname{Pol}(\mathbb{D})_q$;

\item[(iii)] $\pi(\mathsf{f})$ is the identically zero operator on $\operatorname{Pol}(\mathbb{D})_q$;

\item[(iv)] both $\pi(\mathsf{e})$ and $\pi(\mathsf{f})$ are the identically zero operators on $\operatorname{Pol}(\mathbb{D})_q$.
\end{description}
\end{lemma}

{\bf Proof.} Assume (i). Clearly the weight of any monomial in $z$, $z^*$ is $\pm 1$. On the other hand, it follows from \eqref{ke} that $\pi(\mathsf{e})(z)$, if non-zero, should be a weight vector whose weight is $\pm q^2\ne\pm 1$. Hence $\pi(\mathsf{e})(z)=0$. In a similar way, $\pi(\mathsf{e})(z^*)=0$. Thus we conclude that $\pi(\mathsf{e})\equiv 0$, which is just~(ii). The proof of (i) $\Rightarrow$ (iii) is similar.

Assume (ii). An application of \eqref{effe} to $z$ yields $\left(\pi(\mathsf{k})-\pi\left(\mathsf{k}^{-1}\right)\right)(z)=0$, hence $\alpha=\alpha^{-1}$, $\beta=\beta^{-1}$, which is equivalent to (i). The proof of (iii) $\Rightarrow$ (i) is similar, and the rest of implications are clear. \hfill$\blacksquare$

\medskip

The series of symmetries satisfying the equivalent conditions of Lemma \ref{iz} will be called the $\mathbf{(0)}$-series and is described by

\begin{theorem}
The $\mathbf{(0)}$-series consists of the two $U_q(\mathfrak{sl}_2)$-symmetries on $\operatorname{Pol}(\mathbb{D})_q$ given by
\begin{flalign*}
\mathbf{(0+)}: &&\pi(\mathsf{k})(z) &=z,\qquad\pi(\mathsf{k})(z^*)=z^*&
\\ &&\pi(\mathsf{e})(z) &=\pi(\mathsf{e})(z^*)=\pi(\mathsf{f})(z)=\pi(\mathsf{f})(z^*)=0,&
\\
\\ \mathbf{(0-)}: &&\pi(\mathsf{k})(z) &=-z,\qquad\pi(\mathsf{k})(z^*)=-z^*&
\\ &&\pi(\mathsf{e})(z) &=\pi(\mathsf{e})(z^*)=\pi(\mathsf{f})(z)=\pi(\mathsf{f})(z^*)=0,&
\end{flalign*}
which are non-isomorphic.
\end{theorem}

{\bf Proof.} A routine verification establishes that the above formulas extend from the generators of $U_q(\mathfrak{sl}_2)$ and $\operatorname{Pol}(\mathbb{D})_q$ to well-defined symmetries. The symmetries $\mathbf{(0+)}$ and $\mathbf{(0-)}$ are non-isomorphic, because, by Proposition \ref{auts}, any automorphism of $\operatorname{Pol}(\mathbb{D})_q$ commutes with each of the above the actions of $\mathsf{k}$. \hfill$\blacksquare$

\medskip

Let us introduce the notion of {\it grading jump} $\operatorname{GJ}$, to be used to classify the $U_q(\mathfrak{sl}_2)$-symmetries on $\operatorname{Pol}(\mathbb{D})_q$ that break the equivalent properties listed in Lemma \ref{iz}.

\begin{proposition}\label{GJ_d}
Suppose that $\pi$ is a $U_q(\mathfrak{sl}_2)$-symmetry on $\operatorname{Pol}(\mathbb{D})_q$ which does not belong to $(0)$-series. Then there exists a unique non-zero integer $n$ such that for all $k\in\mathbb{Z}$
\begin{equation}\label{jumps}
\pi(\mathsf{e})\mathcal{A}_k\subset\mathcal{A}_{k+n},\qquad
\pi(\mathsf{f})\mathcal{A}_k\subset\mathcal{A}_{k-n}.
\end{equation}
\end{proposition}

{\bf Proof.} First observe that for any non-zero $z^k\varphi(y)\in\mathcal{A}_k$, $k\ge 0$, one has $\mathbf{wt}\left(z^k\varphi(y)\right)=\alpha^k$, and for a non-zero $\psi(y)z^{*k}\in\mathcal{A}_k$, $k\le 0$, one has $\mathbf{wt}\left(\psi(y)z^{*k}\right)=\alpha^k$, with $\alpha$ being the weight constant for $\pi$ as in \eqref{wc}. Since the homogeneous components $\{\mathcal{A}_k\}_{k\in\mathbb{Z}}$ span $\operatorname{Pol}(\mathbb{D})_q$, one deduces that an arbitrary weight vector has weight of the form $\alpha^m$ for some integer $m$.

Another consequence of our assumption on $\pi$ is that $\pi(\mathsf{e})$ is not the identically zero operator. Since $z$, $z^*$ generate $\operatorname{Pol}(\mathbb{D})_q$, either $\pi(\mathsf{e})(z)$ or $\pi(\mathsf{e})(z^*)$ should be non-zero. Let us first assume that $\pi(\mathsf{e})(z)\ne 0$. It follows from \eqref{wc} and \eqref{ke} that $\pi(\mathsf{e})(z)$ is a weight vector whose weight is $q^2\alpha$. On the other hand, by our above observations this weight should be $\alpha^m$ for some integer $m$. Hence with $n=m-1$ one has $\alpha^n=q^2$; in particular, under the assumptions of the Proposition, $\alpha$ should be a root of $q^2$. Clearly $n\ne 0$ since $q$ is not a root of $1$. Thus one deduces that $\alpha$ is also not a root of $1$, together with $q$, and $n$ as above is unique. In particular, the weights of non-zero homogeneous vectors of different degrees are different. Now $\pi(\mathsf{e})\mathcal{A}_k\subset\mathcal{A}_{k+n}$, $k\in\mathbb{Z}$, becomes a consequence of the general form of an element of $\mathcal{A}_k$ and the relation $\alpha^n=q^2$.

Of course, a similar argument also works in the case when $\pi(\mathsf{e})(z^*)\ne 0$. This argument also allows one to derive a unique integer $n$ such that $\alpha^n=q^2$. Even more, if one assumes that both $\pi(\mathsf{e})(z)$ and $\pi(\mathsf{e})(z^*)$ are non-zero, the integer $n$ produced in each of these procedures should be the same, being a unique solution of the same equation $\alpha^n=q^2$.

Now one can reproduce the same argument(s) as above with $\mathsf{e}$ being replaced by $\mathsf{f}$. We get this way $\alpha^n=q^{-2}$, which leads finally to the relation $\pi(\mathsf{f})\mathcal{A}_k\subset\mathcal{A}_{k-n}$. \hfill $\blacksquare$

\begin{definition}\label{GJ_def}
Let $\pi$ be a $U_q(\mathfrak{sl}_2)$-symmetry on $\operatorname{Pol}(\mathbb{D})_q$. If $\pi$ does not belong to the $(0)$-series, we call the (unique) integer $n$ associated to $\pi$ as in Proposition \ref{GJ_d} the {\bf grading jump} ($\operatorname{GJ}$) for $\pi$. In the case when $\pi$ belongs to the $\mathbf{(0)}$-series, we say that $\operatorname{GJ}=0$.
\end{definition}

\begin{proposition}
$\operatorname{GJ}$ is an isomorphism invariant of $U_q(\mathfrak{sl}_2)$-symmetries on $\operatorname{Pol}(\mathbb{D})_q$.
\end{proposition}

{\bf Proof.} Let $\pi$ be a $U_q(\mathfrak{sl}_2)$-symmetry on $\operatorname{Pol}(\mathbb{D})_q$ and $\Psi$ an automorphism of $\operatorname{Pol}(\mathbb{D})_q$ determined by a non-zero constant $\alpha$ as in Proposition \ref{auts}. Clearly $\Psi\left(z^iz^{*j}\right)=\alpha^{i-j}z^iz^{*j}$, hence $\Psi\mathcal{A}_k=\mathcal{A}_k$, $k\in\mathbb{Z}$. This implies that for the isomorphic symmetry $\xi\mapsto\Psi\pi(\xi)\Psi^{-1}$ the relations
$$
\Psi\pi(\mathsf{e})\Psi^{-1}\mathcal{A}_k\subset\mathcal{A}_{k+n},\qquad
\Psi\pi(\mathsf{f})\Psi^{-1}\mathcal{A}_k\subset\mathcal{A}_{k-n}
$$
are equivalent to \eqref{jumps}. \hfill $\blacksquare$

\medskip

Now we are in a position to compute all the $U_q(\mathfrak{sl}_2)$-symmetries on the quantum disc in terms of the grading jump introduced above, such that each value of $\operatorname{GJ}$ labels a series of symmetries, to be denoted as ($\operatorname{GJ}$)-series.

\section{Symmetries with \boldmath$\operatorname{GJ}>0$}\label{GJ+}

Suppose that $\operatorname{GJ}=n>0$ for a symmetry $\pi$, with the weight constant $\alpha$ subject to $\alpha^n=q^2$. In view of \eqref{jumps} we have $\pi(\mathsf{e})(y)=z^np(y)$  for some polynomial $p$. An application of $\pi(\mathsf{e})$ to \eqref{yz} using Definition \ref{symdef}(i), \eqref{def}, \eqref{jumps}, and \eqref{wc} yields
$$
\pi(\mathsf{e})(yz)=
y\pi(\mathsf{e})(z)+\pi(\mathsf{e})(y)\pi(\mathsf{k})(z)=
y\pi(\mathsf{e})(z)+\alpha z^{n+1}p\left(q^{-2}y\right),
$$
$$
\pi(\mathsf{e})\left(q^{-2}zy\right)=
q^{-2}z\pi(\mathsf{e})(y)+q^{-2}\pi(\mathsf{e})(z)\pi(\mathsf{k})(y)=
q^{-2}z^{n+1}p(y)+q^{-2}\pi(\mathsf{e})(z)y.
$$
Since $\pi(\mathsf{e})(z)\in\mathcal{A}_{n+1}$, this implies
$$
q^{-2n-2}\pi(\mathsf{e})(z)y+\alpha z^{n+1}p\left(q^{-2}y\right)=q^{-2}z^{n+1}p(y)+q^{-2}\pi(\mathsf{e})(z)y,
$$
which is equivalent to
\begin{equation}\label{ez0}
\left(q^{-2n}-1\right)\pi(\mathsf{e})(z)y=
z^{n+1}\left[p(y)-\alpha^{n+1}p\left(q^{-2}y\right)\right].
\end{equation}
Since the l.h.s. here is divisible by $y$, we conclude that $p(y)-\alpha^{n+1}p\left(q^{-2}y\right)$ should be also divisible by $y$. With $p(y)=\sum\limits_{i=0}^mp_iy^i$ and $\alpha$ not a root of $1$, the constant term $\left(1-\alpha^{n+1}\right)p_0$ of $p(y)-\alpha^{n+1}p\left(q^{-2}y\right)$ is zero iff $p_0=0$. Thus $p$ is divisible by $y$, so we can now rewrite the expression for $\pi(\mathsf{e})(y)$ in the form
\begin{equation}\label{ey+}
\pi(\mathsf{e})(y)=z^ns_e(y)y,
\end{equation}
with $s_e(y)=\sum\limits_ia_iy^i$ a polynomial. We need also a generalization of \eqref{ey+} as follows.
\begin{multline*}
\pi(\mathsf{e})(y^k)=\sum_{i=0}^{k-1}y^{k-i-1}z^ns_e(y)y^{i+1}=
z^ns_e(y)\left(\sum_{i=0}^{k-1}q^{-2(k-i-1)n}\right)y^k
\\ =q^{-2nk+2n}\sum_{i=0}^{k-1}q^{2ni}z^ns_e(y)y^k=
q^{-2n(k-1)}\frac{1-q^{2nk}}{1-q^{2n}}z^ns_e(y)y^k=
\frac{q^{-2nk}-1}{q^{-2n}-1}z^ns_e(y)y^k,
\end{multline*}
for $k\ge 0$, hence for any polynomial $\varphi$ one has
$$
\pi(\mathsf{e})(\varphi(y))=\left(q^{-2n}-1\right)^{-1}z^ns_e(y)
\left[\varphi\left(q^{-2n}y\right)-\varphi(y)\right].
$$
Furthermore, with $p$ being replaced by $s_e(y)y$, \eqref{ez0} acquires the form
\begin{equation}\label{ez+}
\pi(\mathsf{e})(z)=\left(q^{-2n}-1\right)^{-1}z^{n+1}\left[s_e(y)-\alpha s_e\left(q^{-2}y\right)\right].
\end{equation}
This implies, via a straightforward induction argument, that with $k\ge 0$
$$
\pi(\mathsf{e})\left(z^k\right)=\left(q^{-2n}-1\right)^{-1}z^{n+k}
\left[s_e(y)-\alpha^ks_e\left(q^{-2k}y\right)\right].
$$
Let us apply $\pi(\mathsf{e})$ to \eqref{yz*} using Definition \ref{symdef}(i), \eqref{def}, \eqref{jumps}, \eqref{ey+}, and \eqref{wc}:
$$
\pi(\mathsf{e})(yz^*)=
y\pi(\mathsf{e})(z^*)+\pi(\mathsf{e})(y)\pi(\mathsf{k})(z^*)=
q^{-2(n-1)}\pi(\mathsf{e})(z^*)y+\alpha^{-1}q^2z^{n-1}s_e\left(q^2y\right)
y(1-y),
$$
\begin{multline*}
\pi(\mathsf{e})\left(q^2z^*y\right)=
q^2z^*\pi(\mathsf{e})(y)+q^2\pi(\mathsf{e})(z^*)\pi(\mathsf{k})(y)
\\ =q^2\left(1-q^{-2}y\right)z^{n-1}s_e(y)y+q^2\pi(\mathsf{e})(z^*)y=
q^2 z^{n-1}\left(1-q^{-2n}y\right)s_e(y)y+q^2\pi(\mathsf{e})(z^*)y.
\end{multline*}
This implies
$$
q^{-2(n-1)}\pi(\mathsf{e})(z^*)y+\alpha^{-1}q^2z^{n-1}s_e\left(q^2y\right)
y(1-y)=q^2 z^{n-1}s_e(y)y\left(1-q^{-2n}y\right)+q^2\pi(\mathsf{e})(z^*)y,
$$
which is equivalent to
\begin{equation}\label{ez*+}
\pi(\mathsf{e})(z^*)=\left(q^{-2n}-1\right)^{-1}z^{n-1}
\left[s_e(y)\left(1-q^{-2n}y\right)-\alpha^{-1}s_e\left(q^2y\right)(1-y)\right].
\end{equation}
Again, a straightforward induction argument establishes that with $0<k\le n$
$$
\pi(\mathsf{e})(z^{*k})=\left(q^{-2n}-1\right)^{-1}z^{n-k}
\left[s_e(y)\left(q^{-2n}y;q^2\right)_k-
\alpha^{-k}s_e\left(q^{2k}y\right)\left(y;q^2\right)_k\right].
$$
Here and in what follows, the standard notation
$$(a;q)_n=\prod_{j=0}^{n-1}\left(1-aq^j\right),$$
is used; see, e.g., \cite[p. xiv]{GR}.

Very similar calculations as above can be reproduced for the generator $\mathsf{f}$. We leave routine details to the reader and present here only the outcome.

In view of \eqref{jumps} we have $\pi(\mathsf{f})(y)=r(y)z^{*n}$  for some polynomial $r$. After establishing that $r$ is divisible by $y$, we rewrite this in the form
$$\pi(\mathsf{f})(y)=s_f(y)yz^{*n}$$
for some polynomial $s_f(y)=\sum\limits_ib_iy^i$. Furthermore, with $\varphi$ an arbitrary polynomial
$$
\pi(\mathsf{f})(\varphi(y))=\left(q^{-2n}-1\right)^{-1}s_f(y)
\left[\varphi\left(q^{-2n}y\right)-\varphi(y)\right]z^{*n}.
$$
One also has
\begin{align}
\pi(\mathsf{f})(z) &=\left(q^{-2n}-1\right)^{-1}
\left[s_f(y)\left(1-q^{-2n}y\right)-\alpha^{-1}s_f\left(q^2y\right)(1-y)\right]
z^{*n-1}, \label{fz+}
\\ \pi(\mathsf{f})\left(z^k\right) &=\left(q^{-2n}-1\right)^{-1}
\left[s_f(y)\left(q^{-2n}y;q^2\right)_k-
\alpha^{-k}s_f\left(q^{2k}y\right)\left(y;q^2\right)_k\right]z^{*n-k},\notag
\\ &\;0<k\le n,\label{fzk+}
\\ \pi(\mathsf{f})(z^*) &=\left(q^{-2n}-1\right)^{-1}\left[s_f(y)-\alpha s_f\left(q^{-2}y\right)\right]z^{*n+1}, \label{fz*+}
\\ \pi(\mathsf{f})\left(z^{*k}\right) &=\left(q^{-2n}-1\right)^{-1}
\left[s_f(y)-\alpha^ks_f\left(q^{-2k}y\right)\right]z^{*n+k},\qquad k\ge 0. \notag
\end{align}
Additionally, we will need below an expression for $\pi(\mathsf{f})\left(z^{2n}\right)$, which is formally not covered by \eqref{fzk+}, but is an easy consequence of the latter with $k=n$:
$$
\pi(\mathsf{f})\left(z^{2n}\right)=\left(q^{-2n}-1\right)^{-1}
z^n\left[s_f\left(q^{-2n}y\right)\left(q^{-4n}y;q^2\right)_n-
q^{-4}s_f\left(q^{2n}y\right)\left(y;q^2\right)_n\right].
$$

The above observations allowed us to derive the relations \eqref{ez+}, \eqref{ez*+}, \eqref{fz+}, \eqref{fz*+}, which, together with \eqref{wc} determine (in our present setting $\operatorname{GJ}=n>0$) a symmetry $\pi$ on the distinguished generators of $U_q(\mathfrak{sl}_2)$ and $\operatorname{Pol}(\mathbb{D})_q$ in terms of the parameters $n$, $\alpha$, and the polynomials $s_e$, $s_f$ of one variable. To produce these relations, \eqref{zz*}, \eqref{ke}, \eqref{kf}, \eqref{def}, \eqref{def1} have been used. Certainly, the parameters of a symmetry are not completely arbitrary; in particular, $\alpha^n=q^2$. To adjust finally the parameters and clarify the possible form of the polynomials $s_e$, $s_f$, it is suitable to apply the relation \eqref{effe} to $z^n$. For that, we proceed with computing, using the above formulas. The result of these calculations is formulated as

\begin{theorem}\label{ser_gj+}\hfill
\begin{description}
\item[(i)] There exist no $U_q(\mathfrak{sl}_2)$-symmetries on $\operatorname{Pol}(\mathbb{D})_q$ with $\operatorname{GJ}=n>1$.

\item[(ii)] With $\operatorname{GJ}=1$, there exist two 2-parameter series of $U_q(\mathfrak{sl}_2)$-symmetries on $\operatorname{Pol}(\mathbb{D})_q$ as follows.
\begin{flalign*}
\mathbf{(1a)}: &&\pi(\mathsf{k})(z) &=q^2z, &\pi(\mathsf{k})(z^*) &=q^{-2}z^*,&
\\ &&\pi(\mathsf{e})(y) &=q^{-1}b_0^{-1}zy, &\pi(\mathsf{f})(y) &=\left(b_0y+b_1y^2\right)z^*,&
\\ &&\pi(\mathsf{e})(z) &=qb_0^{-1}z^2, &\pi(\mathsf{f})(z) &=-b_0-b_1y^2,&
\\ &&\pi(\mathsf{e})(z^*) &=-q^{-1}b_0^{-1}, &\pi(\mathsf{f})(z^*) &=q^2b_0z^{*2},&
\\ &&b_0,b_1 &\in\mathbb{C},\qquad b_0\ne 0.&&&
\\
\\ \mathbf{(1b)}: &&\pi(\mathsf{k})(z) &=q^2z, &\pi(\mathsf{k})(z^*) &=q^{-2}z^*,&
\\ &&\pi(\mathsf{e})(y) &=z\left(a_0y+a_1y^2\right), &\pi(\mathsf{f})(y) &=q^{-1}a_0^{-1}yz^*,&
\\ &&\pi(\mathsf{e})(z) &=q^2a_0z^2, &\pi(\mathsf{f})(z) &=-q^{-1}a_0^{-1},&
\\ &&\pi(\mathsf{e})(z^*) &=-a_0-a_1y^2, &\pi(\mathsf{f})(z^*) &=qa_0^{-1}z^{*2},&
\\ &&a_0,a_1 &\in\mathbb{C},\qquad a_0\ne 0.&&&
\end{flalign*}
\end{description}
\end{theorem}

{\bf Proof.}
\begin{equation}\label{efzn+}
\begin{split}
\pi(\mathsf{ef})(z^n) &=\left(q^{-2n}-1\right)^{-1}
\pi(\mathsf{e})\left[s_f(y)\left(q^{-2n}y;q^2\right)_n
-\alpha^{-n}s_f\left(q^{2n}y\right)\left(y;q^2\right)_n\right]
\\ &=\left(q^{-2n}-1\right)^{-2}z^n\left[s_e(y)s_f\left(q^{-2n}y\right)
\left(q^{-4n}y;q^2\right)_n\right.
\\ &\hspace{8.5em}-\left(1+q^{-2}\right)s_e(y)s_f(y)
\left(q^{-2n}y;q^2\right)_n
\\ &\hspace{8.5em}\left.+q^{-2}s_e(y)s_f\left(q^{2n}y\right)
\left(y;q^2\right)_n\right].
\end{split}
\end{equation}
\begin{equation}\label{fezn+}
\begin{split}
\pi(\mathsf{fe})(z^n) &=\left(q^{-2n}-1\right)^{-1}
\pi(\mathsf{f})\left\{z^{2n}\left[s_e(y)-
\alpha^ns_e\left(q^{-2n}y\right)\right]\right\}
\\ &=\left(q^{-2n}-1\right)^{-1}\left\{\pi(\mathsf{f})\left(z^{2n}\right)
\left[s_e(y)-q^2s_e\left(q^{-2n}y\right)\right]\right.
\\ &\hspace{7.5em}\left.+\pi(\mathsf{k})^{-1}\left(z^{2n}\right)
\pi(\mathsf{f})\left[s_e(y)-q^2s_e\left(q^{-2n}y\right)\right]\right\}
\\ &=\left(q^{-2n}-1\right)^{-2}\cdot
\\ &\phantom{=}\cdot\left\{z^n\left[s_e(y)s_f\left(q^{-2n}y\right)
\left(q^{-4n}y;q^2\right)_n\right.\right.
\\ &\hspace{3.5em}-q^{-4}s_e(y)s_f\left(q^{2n}y\right)\left(y;q^2\right)_n
\\ &\hspace{3.5em}-q^2s_e\left(q^{-2n}y\right)s_f\left(q^{-2n}y\right)
\left(q^{-4n}y;q^2\right)_n
\\ &\hspace{3.5em}\left.+q^{-2}s_e\left(q^{-2n}y\right)
s_f\left(q^{2n}y\right)\left(y;q^2\right)_n\right]
\\ &\hspace{1.6em}+\alpha^{-2n}z^{2n}s_f(y)\cdot
\\ &\hspace{2.2em}\left.\cdot\left[s_e\left(q^{-2n}y\right)
-q^2s_e\left(q^{-4n}y\right)-s_e(y)
+q^2s_e\left(q^{-2n}y\right)\right]z^{*n}\right\}
\\ &=\left(q^{-2n}-1\right)^{-2}z^n\left[s_e(y)s_f\left(q^{-2n}y\right)
\left(q^{-4n}y;q^2\right)_n\right.
\\ &\hspace{8.5em}-q^2s_e\left(q^{-2n}y\right)s_f\left(q^{-2n}y\right)
\left(q^{-4n}y;q^2\right)_n
\\ &\hspace{8.5em}-q^{-4}s_e\left(q^{2n}y\right)s_f\left(q^{2n}y\right)
\left(y;q^2\right)_n
\\ &\hspace{8.5em}\left.+q^{-2}s_e(y)s_f\left(q^{2n}y\right)
\left(y;q^2\right)_n\right].
\end{split}
\end{equation}
Finally, we combine \eqref{efzn+} and \eqref{fezn+} to get
\begin{equation}\label{(ef-fe)zn+}
\begin{split}
\pi(\mathsf{ef-fe})(z^n) &=\left(q^{-2n}-1\right)^{-2}z^n\left[
-\left(1+q^{-2}\right)s_e(y)s_f(y)\left(q^{-2n}y;q^2\right)_n\right.
\\ &\hspace{8.5em}+q^2s_e\left(q^{-2n}y\right)s_f\left(q^{-2n}y\right)
\left(q^{-4n}y;q^2\right)_n
\\ &\hspace{8.5em}\left.+q^{-4}s_e\left(q^{2n}y\right)s_f\left(q^{2n}y\right)
\left(y;q^2\right)_n\right].
\end{split}
\end{equation}
In our present context $\operatorname{GJ}>0$ both $s_e$ and $s_f$ are non-zero polynomials, which can be readily deduced (in the case of $s_e$) from \eqref{ez+}, \eqref{ez*+}, and Lemma \ref{iz}; a similar argument also works also in the case of $s_f$. Let $n_e$, $n_f$ be the degrees of $s_e$ and $s_f$, respectively, so that $s_e(y)=a_{n_e}y^{n_e}+\text{(lower terms)}$, $s_f(y)=b_{n_f}y^{n_f}+\text{(lower terms)}$, with $a_{n_e}$, $b_{n_f}$ being non-zero constants.

As a consequence of \eqref{(ef-fe)zn+}, we deduce that $\pi(\mathsf{ef-fe})(z^n)=\left(q^{-2n}-1\right)^{-2}z^nh(y)$, with $h(y)$ being a non-zero polynomial whose highest term is
\begin{multline*}
(-1)^na_{n_e}b_{n_f}\cdot
\\ \cdot\left(-q^{-(n+1)n}-q^{-(n+1)n-2}+q^{-(3n+1)n+2-2n(n_e+n_f)}
+q^{(n-1)n-4+2n(n_e+n_f)}\right)y^{n_e+n_f+n}.
\end{multline*}
This, together with \eqref{effe}, $n_e\ge 0$ $n_f\ge 0$, $n>0$, implies
$$
-q^{-(n+1)n}-q^{-(n+1)n-2}+q^{-(3n+1)n+2-2n(n_e+n_f)}
+q^{(n-1)n-4+2n(n_e+n_f)}=0.
$$
Substituting here $t=q^{2n(n_e+n_f)}$, we get the equation
$$t^2-q^{-2n^2+2}(1+q^2)t+q^{-4n^2+2}=0,$$
whose roots are $t_1=q^{-2n^2+4}$ and $t_2=q^{-2n^2+2}$.

In the first case we have $q^{2n(n_e+n_f)}=q^{-2n^2+4}$, and since $q$ is not a root of $1$, this is equivalent to
\begin{equation}\label{GJ_2+}
n^2+n(n_e+n_f)-2=0.
\end{equation}
This equation with respect to $n$ has 2 real roots, and the conjectured existence of $U_q(\mathfrak{sl}_2)$-symmetries on $\operatorname{Pol}(\mathbb{D})_q$ with some specific values of $n$, $n_e$, $n_f$ should imply that at least one of the two roots $n_1$, $n_2$ is a positive integer. Let it be $n_1$, then $n_2=-2/n_1$ is negative, hence $2/n_1=n_1+n_e+n_f$ is also a positive integer. Assuming in the latter relation $n_1=2$ we get $n_e+n_f=-1$, which is impossible. So it remains the only possibility $n_1=1$, which appears to be a root of \eqref{GJ_2+} iff $n_e+n_f=1$.

A very similar argument establishes that in the second case $q^{2n(n_e+n_f)}=q^{-2n^2+2}$, only the value $n_e+n_f=0$ guarantees the existence of a positive integral root $n$, which is $n=1$.

We conclude that the only positive value of grading jump $\operatorname{GJ}$ under which there exist $U_q(\mathfrak{sl}_2)$-symmetries on $\operatorname{Pol}(\mathbb{D})_q$ is $\operatorname{GJ}=1$; in this case  one should have $\deg s_e+\deg s_f\le 1$.

Our next step is to substitute $n=1$ to \eqref{(ef-fe)zn+} and then to consider the 2 cases as follows: set in \eqref{(ef-fe)zn+} $s_e(y)=a_0$, $s_f(y)=b_0+b_1y$ (respectively, $s_e(y)=a_0+a_1y$, $s_f(y)=b_0$) and apply \eqref{effe} in order to exclude $a_0=q^{-1}b_0^{-1}$ (respectively, $b_0=q^{-1}a_0^{-1}$) in order to obtain finally the series $\mathbf{(1a)}$ (respectively, $\mathbf{(1b)}$), which, already at this point, appear just as in the formulation of the present Theorem. This calculation is completely routine and is left to the reader.

It turns out that one needs not try finding more relations between $a_i$, $b_j$ (e.g., via applying \eqref{effe} to $z^*$). Instead, it suffices to use the formulas for series $\mathbf{(1a)}$ and $\mathbf{(1b)}$ of symmetries as in the formulation of our Theorem in order to apply the generators of $U_q(\mathfrak{sl}_2)$ to the relation in $\operatorname{Pol}(\mathbb{D})_q$, and vice versa, every relation in $U_q(\mathfrak{sl}_2)$ to the generators of $\operatorname{Pol}(\mathbb{D})_q$; in each case one gets the identity. This calculation, while being completely routine (and thus left to the reader), establishes that the formulas for the series $\mathbf{(1a)}$ and $\mathbf{(1b)}$ as in the formulation determine well defined $U_q(\mathfrak{sl}_2)$-symmetries on $\operatorname{Pol}(\mathbb{D})_q$ for all values of the parameters involved therein. \hfill$\blacksquare$

\begin{remark}\label{1is}
The series of symmetries $\mathbf{(1a)}$ and $\mathbf{(1b)}$ are not disjoint. Their intersection is the 1-parameter series determined by setting in $\mathbf{(1a)}$ $b_1=0$; equivalently, it can be produced by setting in $\mathbf{(1b)}$ $a_1=0$ and then substituting $q^{-1}a_0^{-1}=b_0$.
\end{remark}

\section{Symmetries with \boldmath$\operatorname{GJ}<0$}\label{GJ-}

Now assume that $\operatorname{GJ}=-n<0$ for a symmetry $\pi$, with the weight constant $\alpha$ subject to $\alpha^n=q^{-2}$. Although the arguments used below are similar to those applied in Section \ref{GJ+}, one encounters certain diversity in formulas which results in some different conclusions.

In view of \eqref{jumps} we have $\pi(\mathsf{e})(y)=r_e(y)z^{*n}$ and $\pi(\mathsf{f})(y)=z^nr_f(y)$ for some polynomials $r_e$, $r_f$. It turns out that the property of divisibility of $r_e(y)$ and $r_f(y)$ by $y$ does not hold for all $n$ as it was the case in Section \ref{GJ+}. Now let us assume that $n>1$, respectively, $\operatorname{GJ}<-1$. It will be demonstrated below that the above divisibility property should be valid in this case, just as in Section \ref{GJ+}.

An application of $\pi(\mathsf{e})$ to \eqref{yz} using Definition \ref{symdef}(i), \eqref{def}, \eqref{jumps}, and \eqref{wc} yields
$$
\pi(\mathsf{e})(yz)=
y\pi(\mathsf{e})(z)+\pi(\mathsf{e})(y)\pi(\mathsf{k})(z)=
y\pi(\mathsf{e})(z)+\alpha r_e(y)\left(1-q^{-2n}y\right)z^{*n-1},
$$
$$
\pi(\mathsf{e})\left(q^{-2}zy\right)=
q^{-2}z\pi(\mathsf{e})(y)+q^{-2}\pi(\mathsf{e})(z)\pi(\mathsf{k})(y)=
q^{-2}r_e\left(q^2y\right)(1-y)z^{*n-1}+q^{-2n}y\pi(\mathsf{e})(z).
$$
This implies
$$
y\pi(\mathsf{e})(z)+\alpha r_e(y)\left(1-q^{-2n}y\right)z^{*n-1}=
q^{-2}r_e\left(q^2y\right)(1-y)z^{*n-1}+q^{-2n}y\pi(\mathsf{e})(z),
$$
which is equivalent to
\begin{equation}\label{ez1}
\left(q^{-2n}-1\right)y\pi(\mathsf{e})(z)=
\left[\alpha r_e(y)\left(1-q^{-2n}y\right)-q^{-2}r_e\left(q^2y\right)(1-y)\right]z^{*n-1}.
\end{equation}
Since the l.h.s. here is divisible by $y$, we conclude that $\alpha r_e(y)\left(1-q^{-2n}y\right)-q^{-2}r_e\left(q^2y\right)(1-y)$ should be divisible by $y$. With $r_e(y)=\sum\limits_{i=0}^mr_iy^i$ and $\alpha$ not a root of $1$, the constant term $\left(\alpha-q^{-2}\right)r_0$ of $\alpha r_e(y)\left(1-q^{-2n}y\right)-q^{-2}r_e\left(q^2y\right)(1-y)$, under our current assumption $n>1$ is zero iff $r_0=0$. Thus $r_e(y)$ is divisible by $y$, so we can now rewrite the expression for $\pi(\mathsf{e})(y)$ in the form
\begin{equation}\label{ey-}
\pi(\mathsf{e})(y)=s_e(y)yz^{*n},
\end{equation}
with $s_e(y)=\sum\limits_ia_iy^i$ a polynomial. Now \eqref{ey-} can be generalized as follows:
$$
\pi(\mathsf{e})(y^k)=\sum_{i=0}^{k-1}y^{k-i-1}s_e(y)yz^{*n}y^i=
s_e(y)y^k\left(\sum_{i=0}^{k-1}q^{-2in}\right)z^{*n}=
\frac{q^{-2kn}-1}{q^{-2n}-1}s_e(y)y^kz^{*n},
$$
for $k\ge 0$, hence for any polynomial $\varphi$ one has
$$
\pi(\mathsf{e})(\varphi(y))=\left(q^{-2n}-1\right)^{-1}s_e(y)
\left[\varphi\left(q^{-2n}y\right)-\varphi(y)\right]z^{*n}.
$$
Furthermore, with $r_e(y)$ being replaced by $s_e(y)y$, \eqref{ez1} acquires the form
\begin{equation}\label{ez-}
\pi(\mathsf{e})(z)=\left(q^{-2n}-1\right)^{-1}\left[\alpha s_e(y)\left(1-q^{-2n}y\right)-s_e\left(q^2y\right)(1-y)\right]z^{*n-1}.
\end{equation}
This implies, via a straightforward induction argument, that with $0\le k\le n$
\begin{equation}\label{ezk-}
\pi(\mathsf{e})\left(z^k\right)=\left(q^{-2n}-1\right)^{-1}
\left[\alpha^ks_e(y)\left(q^{-2n}y;q^2\right)_k
-s_e\left(q^{2k}y\right)\left(y;q^2\right)_k\right]z^{*n-k}.
\end{equation}
Next, we apply $\pi(\mathsf{e})$ to \eqref{yz*} using Definition \ref{symdef}(i), \eqref{def}, \eqref{jumps}, and \eqref{wc}:
$$
\pi(\mathsf{e})(yz^*)=
y\pi(\mathsf{e})(z^*)+\pi(\mathsf{e})(y)\pi(\mathsf{k})(z^*)=
y\pi(\mathsf{e})(z^*)+\alpha^{-1}s_e(y)yz^{*n+1},
$$
$$
\pi(\mathsf{e})\left(q^2z^*y\right)=
q^2z^*\pi(\mathsf{e})(y)+q^2\pi(\mathsf{e})(z^*)\pi(\mathsf{k})(y)=
s_e\left(q^{-2}y\right)yz^{*n+1}+q^{-2n}y\pi(\mathsf{e})(z^*).
$$
This implies
$$
y\pi(\mathsf{e})(z^*)+\alpha^{-1}s_e(y)yz^{*n+1}=
s_e\left(q^{-2}y\right)yz^{*n+1}+q^{-2n}y\pi(\mathsf{e})(z^*),
$$
which is equivalent to
\begin{equation}\label{ez*-}
\pi(\mathsf{e})(z^*)=\left(q^{-2n}-1\right)^{-1}
\left[\alpha^{-1}s_e(y)-s_e\left(q^{-2}y\right)\right]z^{*n+1}.
\end{equation}
Now an induction argument allows one to establish that with $k\ge 0$
$$
\pi(\mathsf{e})(z^{*k})=\left(q^{-2n}-1\right)^{-1}
\left[\alpha^{-k}s_e(y)-s_e\left(q^{-2k}y\right)\right]z^{*n+k}.
$$

Very similar calculations as above can be reproduced for the generator $\mathsf{f}$. We leave routine details to the reader and present here only the outcome.

In view of \eqref{jumps} we have $\pi(\mathsf{f})(y)=r(y)z^{*n}$  for some polynomial $r$. Again, it turns out that $r(y)$ is divisible by $y$, hence
$$\pi(\mathsf{f})(y)=s_f(y)yz^{*n}$$
for some polynomial $s_f(y)=\sum\limits_ib_iy^i$. Then with $\varphi$ an arbitrary polynomial
\begin{equation}\label{fphi-}
\pi(\mathsf{f})(\varphi(y))=\left(q^{-2n}-1\right)^{-1}z^ns_f(y)
\left[\varphi\left(q^{-2n}y\right)-\varphi(y)\right].
\end{equation}
We compute also
\begin{align}
\pi(\mathsf{f})(z) &=\left(q^{-2n}-1\right)^{-1}z^{n+1}
\left[\alpha^{-1}s_f(y)-s_f\left(q^{-2}y\right)\right], \notag
\\ \pi(\mathsf{f})\left(z^k\right) &=\left(q^{-2n}-1\right)^{-1}z^{n+k}
\left[\alpha^{-k}s_f(y)-s_f\left(q^{-2k}y\right)\right],\qquad k\ge 0,\label{fzk-}
\\ \pi(\mathsf{f})(z^*) &=\left(q^{-2n}-1\right)^{-1}z^{n-1}
\left[\alpha s_f(y)\left(1-q^{-2n}y\right)-
s_f\left(q^2y\right)(1-y)\right], \notag
\\ \pi(\mathsf{f})\left(z^{*k}\right) &=\left(q^{-2n}-1\right)^{-1}z^{n-k}
\left[\alpha^ks_f(y)\left(q^{-2n}y;q^2\right)_k
-s_f\left(q^{2k}y\right)\left(y;q^2\right)_k\right],\notag
\\ &\qquad 0\le k\le n.\notag
\end{align}
We will also need below an expression for $\pi(\mathsf{e})\left(z^{2n}\right)$. It is not covered by \eqref{ezk-}, but easily follows from the latter with $k=n$:
$$
\pi(\mathsf{e})\left(z^{2n}\right)=\left(q^{-2n}-1\right)^{-1}z^n
\left[q^{-4}s_e\left(q^{-2n}y\right)\left(q^{-4n}y;q^2\right)_n
-s_e\left(q^{2n}y\right)\left(y;q^2\right)_n\right].
$$
Similarly to Section \ref{GJ+}, a symmetry $\pi$ in our present setting $\operatorname{GJ}=-n<-1$ (if any) is now determined on the distinguished generators of $U_q(\mathfrak{sl}_2)$ and $\operatorname{Pol}(\mathbb{D})_q$ in terms of the parameters $n$, $\alpha$ ($\alpha^n=q^{-2}$), and the polynomials $s_e$, $s_f$ of one variable. To clarify the very existence of symmetries in this case, it is suitable to apply the relation \eqref{effe} to $z^n$. The outcome is formulated as

\begin{proposition}
There exist no $U_q(\mathfrak{sl}_2)$-symmetries on $\operatorname{Pol}(\mathbb{D})_q$ with $\operatorname{GJ}<-1$.
\end{proposition}
{\bf Proof.}
An application of \eqref{fzk-} (with $k=n$) yields
\begin{equation*}
\begin{split}
&\pi(\mathsf{ef})(z^n)
\\ &=\left(q^{-2n}-1\right)^{-1}
\left\{z^{2n}\pi(\mathsf{e})\left[q^2s_f(y)-
s_f\left(q^{-2n}y\right)\right]\right.
\\ &\hspace{7.5em}\left.+\pi(\mathsf{e})\left(z^{2n}\right)
\pi(\mathsf{k})\left[q^2s_f(y)-s_f\left(q^{-2n}y\right)\right]\right\}
\\ &=\left(q^{-2n}-1\right)^{-2}\cdot
\\ &\quad\cdot\left\{z^{2n}s_e(y)\left[q^2s_f\left(q^{-2n}y\right)
-s_f\left(q^{-4n}y\right)-q^2s_f(y)+s_f\left(q^{-2n}y\right)\right]z^{*n}
\right.
\\ &\qquad\left.+z^n\left[q^{-4}s_e\left(q^{-2n}y\right)
\left(q^{-4n}y;q^2\right)_n-s_e\left(q^{2n}y\right)\left(y;q^2\right)_n
\right]\left[q^2s_f(y)-s_f\left(q^{-2n}y\right)\right]\right\}
\\ &=\left(q^{-2n}-1\right)^{-2}z^n\cdot
\\ &\quad\cdot\left\{s_e\left(q^{2n}y\right)\left[q^2s_f(y)-
s_f\left(q^{-2n}y\right)-q^2s_f\left(q^{2n}y\right)+s_f(y)\right]
\left(y;q^2\right)_n\right.
\\ &\qquad+\left[q^{-2}s_e\left(q^{-2n}y\right)s_f(y)-
q^{-4}s_e\left(q^{-2n}y\right)s_f\left(q^{-2n}y\right)\right]
\left(q^{-4n}y;q^2\right)_n
\\ &\qquad\left.+\left[-q^2s_e\left(q^{2n}y\right)s_f(y)+
s_e\left(q^{2n}y\right)s_f\left(q^{-2n}y\right)\right]
\left(y;q^2\right)_n\right\}
\\ &=\left(q^{-2n}-1\right)^{-2}z^n\cdot
\\ &\quad\cdot\left\{\left[s_e\left(q^{2n}y\right)s_f(y)-
q^2s_e\left(q^{2n}y\right)s_f\left(q^{2n}y\right)\right]\left(y;q^2\right)_n
\right.
\\ &\qquad+\left.\left[q^{-2}s_e\left(q^{-2n}y\right)s_f(y)-
q^{-4}s_e\left(q^{-2n}y\right)s_f\left(q^{-2n}y\right)\right]
\left(q^{-4n}y;q^2\right)_n\right\}.
\end{split}
\end{equation*}
On the other hand, an application of \eqref{ezk-} and \eqref{fphi-} yields

\begin{equation*}
\begin{split}
\pi(\mathsf{fe})(z^n) &=\left(q^{-2n}-1\right)^{-1}
\pi(\mathsf{f})\left[q^{-2}s_e(y)\left(q^{-2n}y;q^2\right)_n
-s_e\left(q^{2n}y\right)\left(y;q^2\right)_n\right]
\\ &=\left(q^{-2n}-1\right)^{-2}z^ns_f(y)
\left[q^{-2}s_e\left(q^{-2n}y\right)\left(q^{-4n}y;q^2\right)_n\right.
\\ &\hspace{10.7em}-s_e(y)\left(q^{-2n}y;q^2\right)_n
\\ &\hspace{10.7em}-q^{-2}s_e(y)\left(q^{-2n}y;q^2\right)_n
\\ &\hspace{10.7em}\left.+s_e\left(q^{2n}y\right)\left(y;q^2\right)_n\right]
\\ &=\left(q^{-2n}-1\right)^{-2}z^n
\left[q^{-2}s_e\left(q^{-2n}y\right)s_f(y)\left(q^{-4n}y;q^2\right)_n\right.
\\ &\hspace{8.5em}-\left(1+q^{-2}\right)s_e(y)s_f(y)\left(q^{-2n}y;q^2\right)_n
\\ &\hspace{8.5em}\left.+s_e\left(q^{2n}y\right)s_f(y)\left(y;q^2\right)_n
\right],
\end{split}
\end{equation*}
whence
\begin{equation}\label{(ef-fe)zn-}
\begin{split}
\pi(\mathsf{ef-fe})(z^n) &=\left(q^{-2n}-1\right)^{-2}z^n\left[-q^2s_e\left(q^{2n}y\right)
s_f\left(q^{2n}y\right)\left(y;q^2\right)_n\right.
\\ &\hspace{8.6em}+\left(1+q^{-2}\right)s_e(y)s_f(y)
\left(q^{-2n}y;q^2\right)_n
\\ &\hspace{8.6em}\left.-q^{-4}s_e\left(q^{-2n}y\right)
s_f\left(q^{-2n}y\right)\left(q^{-4n}y;q^2\right)_n\right].
\end{split}
\end{equation}
In the present case $\operatorname{GJ}<-1$ both $s_e$ and $s_f$ are non-zero polynomials. This can be readily deduced for $s_e$ from \eqref{ez-}, \eqref{ez*-}, and Lemma \ref{iz}; a similar argument works also in the case of $s_f$. Let $n_e$, $n_f$ be the degrees of $s_e$ and $s_f$, respectively, so that $s_e(y)=a_{n_e}y^{n_e}+\text{(lower terms)}$, $s_f(y)=b_{n_f}y^{n_f}+\text{(lower terms)}$, with $a_{n_e}$, $b_{n_f}$ being non-zero constants.

One can observe from \eqref{(ef-fe)zn-} that $\pi(\mathsf{ef-fe})(z^n)=\left(q^{-2n}-1\right)^{-2}z^nh(y)$, where $h(y)$ is a non-zero polynomial whose highest term is
\begin{multline*}
(-1)^na_{n_e}b_{n_f}\cdot
\\ \cdot\left(-q^{2+(n-1)n+2n(n_e+n_f)}+\left(1+q^{-2}\right)q^{(-n-1)n}+
q^{-4+(-3n-1)n-2n(n_e+n_f)}\right)y^{n_e+n_f+n}.
\end{multline*}
This, together with \eqref{effe}, $n_e\ge 0$, $n_f\ge 0$, $n>1$, implies that
$$
-q^{2+(n-1)n+2n(n_e+n_f)}+\left(1+q^{-2}\right)q^{(-n-1)n}+
q^{-4+(-3n-1)n-2n(n_e+n_f)}=0.
$$
Substituting here $t=q^{2n(n_e+n_f)}$, we obtain the equation
$$t^2-(1+q^2)q^{-2n^2-2}t+q^{-4n^2-6}=0,$$
whose roots are $t_1=q^{-2n^2-4}$ and $t_2=q^{-2n^2-2}$.

In the first case we deduce that $n_e$, $n_f$, $n>1$ should be subject to $q^{2n(n_e+n_f)}=q^{-2n^2-4}$, and since $q$ is not a root of $1$, this is equivalent to
$$n^2+(n_e+n_f)n+2=0.$$
Obviously, this equation has no integral solutions $n>1$.

Similarly, we establish in the second case that $q^{2n(n_e+n_f)}=q^{-2n^2-2}$, or, equivalently
$$n^2+(n_e+n_f)n+1=0.$$
Again, this appears to be impossible for integral $n>1$. The Proposition is proved.\hfill$\blacksquare$

\begin{theorem}
With $\operatorname{GJ}=-1$, there exists two 2-parameter series of $U_q(\mathfrak{sl}_2)$-symmetries on $\operatorname{Pol}(\mathbb{D})_q$ as follows.
\begin{flalign*}
\mathbf{(-1a)}:\!\!\!\! &&\pi(\mathsf{k})(z) &=q^{-2}z, &\pi(\mathsf{k})(z^*) &=q^2z^*,&
\\ &&\pi(\mathsf{e})(y) &=-qb_1^{-1}z^*, &\pi(\mathsf{f})(y) &=z(b_0+b_1y),&
\\ &&\pi(\mathsf{e})(z) &=q^{-1}b_1^{-1}, &\pi(\mathsf{f})(z) &=-q^2b_1z^2,&
\\ &&\pi(\mathsf{e})(z^*) &=0, &\pi(\mathsf{f})(z^*) &=-q^{-2}b_0+b_1-\left(1+q^{-2}\right)b_1y,\!\!\!\!&
\\ &&b_0,b_1 &\in\mathbb{C},\qquad b_1\ne 0.&&&
\\
\\ \mathbf{(-1b)}:\!\!\!\! &&\pi(\mathsf{k})(z) &=q^{-2}z, &\pi(\mathsf{k})(z^*) &=q^2z^*,&
\\ &&\pi(\mathsf{e})(y) &=(a_0+a_1y)z^*, &\pi(\mathsf{f})(y) &=-qa_1^{-1}z,&
\\ &&\pi(\mathsf{e})(z) &=-q^{-2}a_0+a_1-\left(1+q^{-2}\right)a_1y,\!\!\! &\pi(\mathsf{f})(z) &=0,&
\\ &&\pi(\mathsf{e})(z^*) &=-q^2a_1z^{*2}, &\pi(\mathsf{f})(z^*) &=q^{-1}a_1^{-1},&
\\ &&a_0,a_1 &\in\mathbb{C},\qquad a_1\ne 0.&&&
\end{flalign*}
\end{theorem}

{\bf Proof.} In this case, with $\pi$ being a   $U_q(\mathfrak{sl}_2)$-symmetry on $\operatorname{Pol}(\mathbb{D})_q$ (if exists), we have $\pi(\mathsf{k})(z)=q^{-2}z$, $\pi(\mathsf{k})(z^*)=q^2z^*$, $\pi(\mathsf{e})(y)=r_e(y)z^*$, $\pi(\mathsf{f})(y)=zr_f(y)$ in view of Definition \ref{GJ_def} and Proposition \ref{GJ_d}. Unlike the cases considered before, we can not claim now that $r_e(y)$ and/or $r_f(y)$ is divisible by $y$. On the other hand, we need to deal with polynomials which are divisible; for that, we introduce the division map $\tau:\mathbb{C}[y]y\to\mathbb{C}[y]$, $\tau[y\mapsto\varphi(y)y](y)=\varphi(y)$. Since everything is embedded into $\operatorname{Pol}(\mathbb{D})_q$, the map $\tau$ is well defined, because $\operatorname{Pol}(\mathbb{D})_q$ is a domain. The following completely obvious property of $\tau$ is going to be useful in what follows:
\begin{equation}\label{div}
\tau(\psi)\circ\beta=\beta^{-1}\tau(\psi\circ\beta),\qquad
\psi\in\mathbb{C}[y]y,\;\beta\in\mathbb{C},
\end{equation}
where $\varphi\circ\beta(y)=\varphi(\beta y)$.

Firstly, one uses a straightforward induction argument in order to establish that for any polynomial $\varphi$ one has
\begin{equation}\label{ephi-1}
\pi(\mathsf{e})(\varphi(y))=\left(q^{-2}-1\right)^{-1}r_e(y)
\tau\left[\varphi\left(q^{-2}y\right)-\varphi(y)\right]z^*.
\end{equation}

Next, in order to produce the expression for $\pi(\mathsf{e})(z)$, we compute
$$
\pi(\mathsf{e})(yz)=
y\pi(\mathsf{e})(z)+\pi(\mathsf{e})(y)\pi(\mathsf{k})(z)=
y\pi(\mathsf{e})(z)+q^{-2}r_e(y)\left(1-q^{-2}y\right),
$$
$$
\pi(\mathsf{e})\left(q^{-2}zy\right)=
q^{-2}z\pi(\mathsf{e})(y)+q^{-2}\pi(\mathsf{e})(z)\pi(\mathsf{k})(y)=
q^{-2}r_e\left(q^2y\right)(1-y)+q^{-2}y\pi(\mathsf{e})(z).
$$
This implies
$$
y\pi(\mathsf{e})(z)+q^{-2}r_e(y)\left(1-q^{-2}y\right)z^{*n-1}=
q^{-2}r_e\left(q^2y\right)(1-y)+q^{-2}y\pi(\mathsf{e})(z),
$$
which is equivalent to
\begin{equation}\label{ez-1p}
\left(q^{-2}-1\right)y\pi(\mathsf{e})(z)=
\left[q^{-2}r_e(y)\left(1-q^{-2}y\right)-q^{-2}r_e\left(q^2y\right)(1-y)\right].
\end{equation}
It is easy to observe that the constant term of the polynomial in the r.h.s of \eqref{ez-1p} is zero, so the polynomial is divisible by $y$, whence
\begin{align}
\pi(\mathsf{e})(z) &=\left(q^{-2}-1\right)^{-1}q^{-2}\tau\left[ r_e(y)\left(1-q^{-2}y\right)-r_e\left(q^2y\right)(1-y)\right],&
\label{ez-1}
\\ \pi(\mathsf{e})(z^2) &=\left(q^{-2}-1\right)^{-1}q^{-2}z
\tau\left[r_e\left(q^{-2}y\right)\left(1-q^{-4}y\right)-
r_e\left(q^2y\right)(1-y)\right].&\label{ez2-1}
\end{align}
In a similar way, we compute:
$$
\pi(\mathsf{e})(yz^*)=
y\pi(\mathsf{e})(z^*)+\pi(\mathsf{e})(y)\pi(\mathsf{k})(z^*)=
y\pi(\mathsf{e})(z^*)+q^2r_e(y)z^{*2},
$$
$$
\pi(\mathsf{e})\left(q^2z^*y\right)=
q^2z^*\pi(\mathsf{e})(y)+q^2\pi(\mathsf{e})(z^*)\pi(\mathsf{k})(y)=
q^2r_e\left(q^{-2}y\right)z^{*2}+q^{-2}y\pi(\mathsf{e})(z^*).
$$
This implies
$$
y\pi(\mathsf{e})(z^*)+q^2r_e(y)yz^{*2}=
q^2r_e\left(q^{-2}y\right)z^{*2}+q^{-2}y\pi(\mathsf{e})(z^*),
$$
which, in view of divisibility of $q^2r_e(y)-q^2r_e\left(q^{-2}y\right)$ by $y$, is equivalent to
\begin{equation}\label{ez*-1}
\pi(\mathsf{e})(z^*)=\left(q^{-2}-1\right)^{-1}q^2
\tau\left[r_e(y)-r_e\left(q^{-2}y\right)\right]z^{*2}.
\end{equation}

The above calculations can be reproduced for the generator $\mathsf{f}$. The details are left to the reader; the outcome only is given below.
\begin{align}
\pi(\mathsf{f})(\varphi(y)) &=\left(q^{-2}-1\right)^{-1}zr_f(y)
\tau\left[\varphi\left(q^{-2}y\right)-\varphi(y)\right],\qquad
\varphi\in\mathbb{C}[y],&\label{fphi-1}
\\ \pi(\mathsf{f})(z) &=\left(q^{-2}-1\right)^{-1}q^2z^2
\tau\left[r_f(y)-r_f\left(q^{-2}y\right)\right],&\label{fz-1}
\\ \pi(\mathsf{f})(z^*) &=\left(q^{-2}-1\right)^{-1}q^{-2}
\tau\left[r_f(y)\left(1-q^{-2}y\right)-r_f\left(q^2y\right)(1-y)\right].&
\label{fz*-1}
\end{align}

What remains is to compute the general form of the polynomials $r_e$, $r_f$. To do that, we apply \eqref{effe} to $z$, using \eqref{div} -- \eqref{fz*-1}.
\begin{equation*}
\begin{split}
&\pi(\mathsf{ef})(z)
\\ &=\left(q^{-2}-1\right)^{-1}
\pi(\mathsf{e})\left\{z^2q^2\tau\left[r_f(y)-
r_f\left(q^{-2}y\right)\right]\right\}
\\ &=\left(q^{-2}-1\right)^{-1}
\left\{z^2q^2\pi(\mathsf{e})\left[\tau\left(r_f(y)-
r_f\left(q^{-2}y\right)\right)\right]\right.
\\ &\left.\hspace{7em}+\pi(\mathsf{e})\left(z^2\right)\pi(\mathsf{k})
\tau\left[q^2r_f(y)-q^2r_f\left(q^{-2}y\right)\right]\right\}
\\ &=\left(q^{-2}-1\right)^{-2}\left\{z^2q^2r_e(y)
\left(q^2\tau\left[\tau\left(q^2r_f\left(q^{-2}y\right)
-q^2r_f\left(q^{-4}y\right)\right)-
\tau\left(r_f(y)-r_f\left(q^{-2}y\right)\right)\right]z^*\right.\right.
\\ &\hspace{7em}\left.+z\tau\left[r_e\left(q^{-2}y\right)
\left(1-q^{-4}y\right)-r_e\left(q^2y\right)(1-y)\right]\cdot
\tau\left[r_f(y)-r_f\left(q^{-2}y\right)\right]\right\}
\\ &=\left(q^{-2}-1\right)^{-2}z
\left\{r_e\left(q^2y\right)\tau^2\left[r_f(y)-r_f\left(q^{-2}y\right)-
q^{-2}r_f\left(q^2y\right)+q^{-2}r_f(y)\right](1-y)\right.
\\ &\left.\hspace{7.6em}+\tau^2\left[\left(r_e\left(q^{-2}y\right)
\left(1-q^{-4}y\right)-r_e\left(q^2y\right)(1-y)\right)\cdot
\left(r_f(y)-r_f\left(q^{-2}y\right)\right)\right]\right\}
\\ &=\left(q^{-2}-1\right)^{-2}z\tau^2
\left[r_e\left(q^2y\right)r_f\left(q^2y\right)(-q^{-2}+q^{-2}y)+
r_e\left(q^2y\right)r_f(y)(q^{-2}-q^{-2}y)\right.
\\ &\left.\hspace{8.5em}+r_e\left(q^{-2}y\right)r_f(y)\left(1-q^{-4}y\right)+
r_e\left(q^{-2}y\right)r_f\left(q^{-2}y\right)\left(-1+q^{-4}y\right)\right].
\end{split}
\end{equation*}
On the other hand, an application of \eqref{ez-1} and \eqref{fphi-1} yields
\begin{equation*}
\begin{split}
\pi(\mathsf{fe})(z) &=\left(q^{-2}-1\right)^{-1}q^{-2}
\pi(\mathsf{f})\left\{\tau\left[r_e(y)\left(1-q^{-2}y\right)
-r_e\left(q^2y\right)(1-y)\right]\right\}
\\ &=\left(q^{-2}-1\right)^{-2}q^{-2}zr_f(y)
\tau\left\{q^2\tau\left[r_e\left(q^{-2}y\right)\left(1-q^{-4}y\right)
-r_e(y)\left(1-q^{-2}y\right)\right]\right.
\\ &\hspace{10.7em}\left.-\tau\left[r_e(y)\left(1-q^{-2}y\right)
-r_e\left(q^{2}y\right)(1-y)\right]\right\}
\\ &=\left(q^{-2n}-1\right)^{-2}z\tau^2
\left[r_e\left(q^{-2}y\right)r_f(y)\left(1-q^{-4}y\right)\right.
\\ &\hspace{9em}+r_e(y)r_f(y)\left(-1-q^{-2}+\left(q^{-2}+q^{-4}\right)y\right)
\\ &\hspace{9em}\left.+r_e\left(q^2y\right)r_f(y)\left(q^{-2}-q^{-2}y\right)\right],
\end{split}
\end{equation*}
whence
\begin{equation}\label{(ef-fe)z-1}
\begin{split}
\pi(\mathsf{ef-fe})(z) &=\left(q^{-2}-1\right)^{-2}z\tau^2
\left[r_e\left(q^2y\right)r_f\left(q^2y\right)
\left(-q^{-2}+q^{-2}y\right)\right.
\\ &\hspace{8.6em}+r_e\left(q^{-2}y\right)r_f\left(q^{-2}y\right)
\left(-1+q^{-4}y\right)
\\ &\hspace{8.6em}\left.+r_e(y)r_f(y)
\left(1+q^{-2}-\left(q^{-2}+q^{-4}\right)y\right)\right].
\end{split}
\end{equation}

In the present case $\operatorname{GJ}=-1$ both $r_e$ and $r_f$ are non-zero polynomials. One can deduce this for $r_e$ from \eqref{ez-1}, \eqref{ez*-1}, and Lemma \ref{iz}; a similar argument works also in the case of $r_f$. Let $n_e$, $n_f$ be the degrees of $r_e$ and $r_f$, respectively, so that $r_e(y)=a_{n_e}y^{n_e}+\text{(lower terms)}$, $r_f(y)=b_{n_f}y^{n_f}+\text{(lower terms)}$, with $a_{n_e}$, $b_{n_f}$ being non-zero constants.

Let us rewrite, in view of the divisibility issues described above, \eqref{(ef-fe)z-1} in the form $\pi(\mathsf{ef-fe})(z)=\left(q^{-2n}-1\right)^{-2}zh(y)$, where $h(y)$ is a non-zero polynomial. The highest term of $h(y)$ is
$$
a_{n_e}b_{n_f}\left(-q^{2(n_e+n_f)-2}+q^{-2(n_e+n_f)-4}-q^{-2}-q^{-4}\right)
y^{n_e+n_f-1}.
$$
Of course, one has here $n_e+n_f\ge 1$, because otherwise ($n_e=n_f=0$) one deduces from \eqref{(ef-fe)z-1} that $\pi(\mathsf{ef-fe})(z)=0$, which, in view of \eqref{effe}, implies $\alpha=\pm 1$, contradicting $\operatorname{GJ}=-1$.

Assuming $n_e+n_f>1$, one clearly observes from \eqref{effe} applied to $z$ that
$$-q^{2(n_e+n_f)-2}+q^{-2(n_e+n_f)-4}-q^{-2}-q^{-4}=0.$$
To find the possible values of $n_e+n_f>1$ that could make possible the latter relation, we substitute here $t=q^{2n(n_e+n_f)}$ in order to get the equation
$$t^2-\left(1+q^{-2}\right)t+q^{-2}=0,$$
whose roots are $t_1=q^{-2}$ and $t_2=1$. Respectively, this yields $n_e+n_f=-1\text{\ or\ }0$, breaking the assumption $n_e+n_f>1$.

Thus we conclude that the only possibility is $n_e+n_f=\deg r_e+\deg r_f=1$.

The final step in producing the series $\mathbf{(-1a)}$ (respectively, $\mathbf{(-1b)}$) as in the formulation of Theorem, is to consider the 2 cases as follows: set in \eqref{(ef-fe)z-1} $r_e(y)=a_0$, $r_f(y)=b_0+b_1y$ (respectively, $r_e(y)=a_0+a_1y$, $r_f(y)=b_0$) and apply \eqref{effe} to exclude $a_0=-qb_1^{-1}$ (respectively, $b_0=-qa_1^{-1}$). This calculation is completely routine and is left to the reader.

Now we suggest to reproduce the final step of the proof of Theorem \ref{ser_gj+}. Namely, instead of searching for more relations between $a_i$, $b_j$ (e.g., via applying \eqref{effe} to $z^*$), it suffices to use the formulas for series $\mathbf{(-1a)}$ and $\mathbf{(-1b)}$ of symmetries as in the formulation of the present Theorem in order to apply the generators of $U_q(\mathfrak{sl}_2)$ to the relation in $\operatorname{Pol}(\mathbb{D})_q$, and vice versa, every relation in $U_q(\mathfrak{sl}_2)$ to the generators of $\operatorname{Pol}(\mathbb{D})_q$. In all the cases one gets the identity, which demonstrates that the formulas for the series $\mathbf{(-1a)}$ and $\mathbf{(-1b)}$ as in the formulation determine well defined $U_q(\mathfrak{sl}_2)$-symmetries on $\operatorname{Pol}(\mathbb{D})_q$ for all values of the parameters involved therein. Again, this calculation appears to be purely technical and thus left to the reader. \hfill$\blacksquare$

\begin{remark}
No symmetry of the series $\mathbf{(-1a)}$ or $\mathbf{(-1b)}$ is isomorphic to a symmetry of the series $\mathbf{(1a)}$ $\mathbf{(1b)}$. This is due to the fact that an arbitrary automorphism of $\operatorname{Pol}(\mathbb{D})_q$ (see Proposition \ref{auts}) commutes with the action of $\pi(\mathsf{k})$ for any $U_q(\mathfrak{sl}_2)$-symmetry $\pi$ on $\operatorname{Pol}(\mathbb{D})_q$.
\end{remark}

\begin{remark}
The series of symmetries $\mathbf{(-1a)}$ and $\mathbf{(-1b)}$ are disjoint. To see this, one can, e.g., observe that in the series $\mathbf{(-1a)}$, $\pi(\mathsf{f})(z)$ is non-zero for any (non-zero) value of the parameter $b_1$; on the other hand, in the series $\mathbf{(-1b)}$ one has $\pi(\mathsf{f})(z)=0$  for all admissible values of the parameters.

A very similar argument can be also used to establish that no symmetry of the series $\mathbf{(-1a)}$ is isomorphic to a symmetry of the series $\mathbf{(-1b)}$. In fact, with $\pi$ being a $\mathbf{(-1b)}$-symmetry and $\Psi$ an automorphism of $\operatorname{Pol}(\mathbb{D})_q$, one readily computes, using Proposition \ref{auts}, that $\Psi\pi(\mathsf{f})\Psi^{-1}(z)=0$.
\end{remark}

\section{A note on involutions}\label{invs}

The approach used above that ignored the presence of involutions both in $\operatorname{Pol}(\mathbb{D})_q$ and in $U_q(\mathfrak{sl}_2)$ was helpful in describing the utmost collection of $U_q(\mathfrak{sl}_2)$-symmetries on $\operatorname{Pol}(\mathbb{D})_q$. However, it would be unnatural to avoid even mentioning at least the straightforward involution on $\operatorname{Pol}(\mathbb{D})_q$, which sends $z$ to $z^*$. As for $U_q(\mathfrak{sl}_2)$, the picture is less plausible, as the latter Hopf algebra admits several involutions (real forms) compatible with the structures on $U_q(\mathfrak{sl}_2)$ as a Hopf algebra. Also, a sort of compatibility is assumed implicit on involutions involved for a $U_q(\mathfrak{sl}_2)$-symmetry on $\operatorname{Pol}(\mathbb{D})_q$. So, we start with recalling relevant definitions, see, e.g., \cite{KS}.

Let $H$ be a Hopf algebra whose comultiplication is $\Delta$, counit is $\varepsilon$, and antipode is $\mathsf{S}$. Suppose $H$ is equipped with an involution $*$, which is an antilinear antiisomorphism. $H$ is called a Hopf $*$-algebra if the following conditions are satisfied. $\Delta:H\to H\otimes H$ is a $*$-homomorphism. The latter means that $\Delta(a^*)=\Delta(a)^*$ for $a\in H$, where the involution of $H\otimes H$ is defined by $(a\otimes b)^*=a^*\otimes b^*$. This definition already implies certain relations between $*$, $\mathsf{S}$, and $\varepsilon$ \cite[1.2.7]{KS}.

Now let $A$ be a unital involutive algebra, whose unit is $\mathbf{1}$, and the involution is denoted by the same symbol $*$ as above. Let also $\pi$ be an $H$-symmetry on $A$. In this specific case, the following compatibility assumption on involutions is implicit \cite{V} (see also \cite{SSV1, SSV2, SZ}):
\begin{equation}\label{invcomp}
(\pi(\xi)a)^*=\pi(\mathsf{S}(\xi)^*)a^*,\qquad\xi\in H,\quad a\in A.
\end{equation}
Here the symmetry sign $\pi$ is used explicitly, unlike \cite{V} where the symmetry in this compatibility property is implicit and thus omitted. Just as in the context of \cite{V}, this part of the definition of a symmetry (structure of $H$-module algebra on $A$) allows a proper application of the involution(s) to the relation
$$\pi(\xi\eta)a=\pi(\xi)(\pi(\eta)a),\qquad\xi,\eta\in H,\quad a\in A.$$

In our specific case $H=U_q(\mathfrak{sl}_2)$, $A=\operatorname{Pol}(\mathbb{D})_q$, we already have a complete list of $U_q(\mathfrak{sl}_2)$-symmetries on $\operatorname{Pol}(\mathbb{D})_q$ described in Sections \ref{trl}, \ref{GJ+}, \ref{GJ-}. What remains  is to extract a (sub)list of symmetries compatible with involutions as in \eqref{invcomp}.

We restrict our considerations to the involution $z\mapsto z^*$ in $\operatorname{Pol}(\mathbb{D})_q$, and reproduce below the list of involutions that make $U_q(\mathfrak{sl}_2)$ a Hopf $*$-algebra \cite[3.1.4]{KS}. The list is exhaustive and contains representatives of equivalence classes (of involutions that can be intertwined by automorphisms of the Hopf algebra $U_q(\mathfrak{sl}_2)$). Each item in this list is related to a specific set of values for $q$, and we additionally keep our initial assumption that $q$ is not a root of $1$.
\begin{description}
\item[(A)] This involution is valid with $q\in\mathbb{R}$, and the corresponding Hopf $*$-algebra $U_q(\mathfrak{su}_2)$ is called the compact real form of $U_q(\mathfrak{sl}_2)$. Explicitly,
$$
\mathsf{k}^*=\mathsf{k},\qquad\mathsf{e}^*=\mathsf{f}\mathsf{k},\qquad
\mathsf{f}^*=\mathsf{k}^{-1}\mathsf{e}.
$$

\item[(B)] Similarly to the previous case, $q\in\mathbb{R}$, and the corresponding Hopf $*$-algebra is denoted by $U_q(\mathfrak{su}_{1,1})$. Explicitly,
$$
\mathsf{k}^*=\mathsf{k},\qquad\mathsf{e}^*=-\mathsf{f}\mathsf{k},\qquad
\mathsf{f}^*=-\mathsf{k}^{-1}\mathsf{e}.
$$

\item[(C)] Let $|q|=1$. The single equivalence class of involutions that make $U_q(\mathfrak{sl}_2)$ a Hopf $*$-algebra is represented by
$$
\mathsf{k}^*=\mathsf{k},\qquad\mathsf{e}^*=\mathsf{e},\qquad
\mathsf{f}^*=\mathsf{f}.
$$
This real form is denoted by $U_q(\mathfrak{sl}_2(\mathbb{R}))$.

\item[(D)] Let $q\in i\mathbb{R}$. An equivalence class of involutions that have no classical counterpart is represented by
$$
\mathsf{k}^*=\mathsf{k},\qquad\mathsf{e}^*=i\mathsf{f}\mathsf{k},\qquad
\mathsf{f}^*=i\mathsf{k}^{-1}\mathsf{e}.
$$

\item[(E)] Again with $q\in i\mathbb{R}$, there exists just one more equivalence class of involutions that have no classical counterpart; it is represented by
$$
\mathsf{k}^*=\mathsf{k},\qquad\mathsf{e}^*=-i\mathsf{f}\mathsf{k},\qquad
\mathsf{f}^*=-i\mathsf{k}^{-1}\mathsf{e}.
$$
\end{description}

Now we are in a position to produce a complete list of series of $U_q(\mathfrak{sl}_2)$-symmetries on $\operatorname{Pol}(\mathbb{D})_q$ which, under presence of involutions both in $U_q(\mathfrak{sl}_2)$ and in $\operatorname{Pol}(\mathbb{D})_q$, admit the compatibility condition \eqref{invcomp}. We start with the following Lemma which describes the cases when \eqref{invcomp} agrees with the algebraic structures both on $U_q(\mathfrak{sl}_2)$ and on $\operatorname{Pol}(\mathbb{D})_q$, regardless of the explicit form of symmetries.

\begin{lemma}\label{icex}
Suppose that the involution on $U_q(\mathfrak{sl}_2)$ is $\mathbf{(A)}$, $\mathbf{(B)}$, $\mathbf{(D)}$, or $\mathbf{(E)}$. Assume also that the relation \eqref{invcomp} is true with $\xi$ taking values in the distinguished set of generators $\mathsf{k},\mathsf{k}^{-1},\mathsf{e},\mathsf{f}\in U_q(\mathfrak{sl}_2)$ and, $a$ being $z$ or $z^*\in\operatorname{Pol}(\mathbb{D})_q$. Then \eqref{invcomp} is true for arbitrary $\xi\in U_q(\mathfrak{sl}_2)$, $a\in\operatorname{Pol}(\mathbb{D})_q$.
\end{lemma}

{\bf Proof.} This Lemma does not allude to an explicit form of a symmetry $\pi$ as in \eqref{invcomp}, so we omit the very symbol $\pi$ throughout the present proof, thus making a symmetry implicit.

Let $a,b=z\text{\ or\ }z^*\in\operatorname{Pol}(\mathbb{D})_q$. Let us now restrict our considerations to the involution $\mathbf{(A)}$ on $U_q(\mathfrak{sl}_2)$. In this case
\begin{align*}
\mathsf{S}(\mathsf{k})^* &=(\mathsf{k}^{-1})^*=\mathsf{k}^{-1},\qquad
\mathsf{S}(\mathsf{k}^{-1})^*=(\mathsf{k})^*=\mathsf{k},&
\\ \mathsf{S}(\mathsf{e})^* &=(-\mathsf{e}\mathsf{k}^{-1})^*=
-\mathsf{k}^{-1}\mathsf{f}\mathsf{k}=-q^2\mathsf{f},&
\\ \mathsf{S}(\mathsf{f})^* &=(-\mathsf{k}\mathsf{f})^*=
-\mathsf{k}^{-1}\mathsf{e}\mathsf{k}=-q^{-2}\mathsf{e}.&
\end{align*}
With this, under the assumptions of Lemma one has
\begin{multline*}
\mathsf{k}(ab)^*=(\mathsf{k}(a)\mathsf{k}(b))^*=
\mathsf{k}(b)^*\mathsf{k}(a)^*=
(\mathsf{S}(\mathsf{k})^*b^*)(\mathsf{S}(\mathsf{k})^*a^*)=
\mathsf{k}^{-1}(b^*)\mathsf{k}^{-1}(a^*)=\mathsf{k}^{-1}(b^*a^*)
\\ =\mathsf{S}(\mathsf{k})^*(ab)^*,
\end{multline*}
and, in a similar way
$$\mathsf{k}^{-1}(ab)^*=\mathsf{S}(\mathsf{k}^{-1})^*(ab)^*.$$
Also, we compute
\begin{multline*}
\mathsf{e}(ab)^*=(a\mathsf{e}(b))^*+(\mathsf{e}(a)\mathsf{k}(b))^*=
\mathsf{e}(b)^*a^*+\mathsf{k}(b)^*\mathsf{e}(a)^*
\\ =(\mathsf{S}(\mathsf{e})^*b^*)a^*+
(\mathsf{S}(\mathsf{k})^*b^*)(\mathsf{S}(\mathsf{e})^*a^*)=
-q^2\mathsf{f}(b^*)a^*-\mathsf{k}^{-1}(b^*)q^2\mathsf{f}(a^*)
\\ =-q^2\Delta(\mathsf{f})(b^*\otimes a^*)=
-q^2\mathsf{f}(b^*a^*)=\mathsf{S}(\mathsf{e})^*(ab)^*,
\end{multline*}
\begin{multline*}
\mathsf{f}(ab)^*=(\mathsf{f}(a)b)^*+(\mathsf{k}^{-1}(a)\mathsf{f}(b))^*=
b^*\mathsf{f}(a^*)+\mathsf{f}(b)^*\mathsf{k}^{-1}(a)^*
\\ =b^*(\mathsf{S}(\mathsf{f})^*a^*)+
(\mathsf{S}(\mathsf{f})^*b^*)(\mathsf{S}(\mathsf{k}^{-1})^*a^*)=
-b^*q^{-2}\mathsf{e}(a^*)-q^{-2}\mathsf{e}(b^*)\mathsf{k}(a^*)
\\ =-q^{-2}\Delta(\mathsf{e})(b^*\otimes a^*)=
-q^2\mathsf{e}(b^*a^*)=\mathsf{S}(\mathsf{f})^*(ab)^*,
\end{multline*}

Similar arguments work in the cases of the rest of involutions on $U_q(\mathfrak{sl}_2)$ listed in the formulation of Lemma. This proves \eqref{invcomp} for $\xi=\mathsf{k},\mathsf{k}^{-1},\mathsf{e},\text{\ or\ }\mathsf{f}\in U_q(\mathfrak{sl}_2)$, $a\in\operatorname{Pol}(\mathbb{D})_q$, due to the anti-linearity in $a$ of the l.h.s. and the r.h.s. of \eqref{invcomp} with a fixed $\xi\in U_q(\mathfrak{sl}_2)$.

Now let $\xi,\eta=\mathsf{k},\mathsf{k}^{-1},\mathsf{e},\text{\ or\ }\mathsf{f}\in U_q(\mathfrak{sl}_2)$, $a\in\operatorname{Pol}(\mathbb{D})_q$. In view of the above observations,
$$
((\xi\eta)a)^*=(\xi(\eta a))^*=\mathsf{S}(\xi)^*(\eta a)^*=\mathsf{S}(\xi)^*\mathsf{S}(\eta)^*a^*=
(\mathsf{S}(\eta)\mathsf{S}(\xi))^*a^*=\mathsf{S}(\xi\eta)^*a^*,
$$
which finishes the proof, due to the anti-linearity in $\xi$ of the l.h.s. and the r.h.s. of \eqref{invcomp} with a fixed $a\in\operatorname{Pol}(\mathbb{D})_q$.\hfill$\blacksquare$

\begin{remark}
The involution $\mathbf{(C)}$ is not covered by Lemma \ref{icex}, because under this involution \eqref{invcomp} fails even on the generators, unless the symmetry in question is either $\mathbf{(0+)}$ or $\mathbf{(0-)}$.
\end{remark}

\begin{theorem}\hfill
\begin{description}
\item[(i)] The symmetries $\mathbf{(0+)}$ and $\mathbf{(0-)}$ admit compatibility for each of the involutions $\mathbf{(A)}$ -- $\mathbf{(E)}$ with the involution in $\operatorname{Pol}(\mathbb{D})_q$.

\item[(ii)] Suppose that $q<0$. Then the involution $\mathbf{(A)}$ possesses the compatibility property \eqref{invcomp} with the involution in $\operatorname{Pol}(\mathbb{D})_q$ so that the latter admits $U_q(\mathfrak{su}_2)$-symmetries under a part of the series $\mathbf{(1a)}$ distinguished by setting there $b_1=0$, $|b_0|^2=-q^{-3}$. Explicitly, those symmetries are
\begin{align*}
\pi(\mathsf{k})(z) &=q^2z, &\pi(\mathsf{k})(z^*) &=q^{-2}z^*,&
\\ \pi(\mathsf{e})(y) &=q^{-1}b_0^{-1}zy, &\pi(\mathsf{f})(y) &=b_0yz^*,&
\\ \pi(\mathsf{e})(z) &=qb_0^{-1}z^2, &\pi(\mathsf{f})(z) &=-b_0,&
\\ \pi(\mathsf{e})(z^*) &=-q^{-1}b_0^{-1}, &\pi(\mathsf{f})(z^*) &=q^2b_0z^{*2},&
\\ b_0 &\in\mathbb{C},\quad |b_0|^2=-q^{-3}.&&&
\end{align*}
In view of Remark \ref{1is}, this set of symmetries is also a part of the series $\mathbf{(1b)}$ distinguished by setting there $a_1=0$, $|a_0|^2=-q$.

\item[(iii)] Suppose that $q>0$. Then the involution $\mathbf{(B)}$ possesses the compatibility property \eqref{invcomp} with the involution in $\operatorname{Pol}(\mathbb{D})_q$ so that the latter admits $U_q(\mathfrak{su}_{1,1})$-symmetries under a part of the series $\mathbf{(1a)}$ distinguished by setting there $b_1=0$, $|b_0|^2=q^{-3}$. Explicitly, those symmetries are
\begin{align*}
\pi(\mathsf{k})(z) &=q^2z, &\pi(\mathsf{k})(z^*) &=q^{-2}z^*,&
\\ \pi(\mathsf{e})(y) &=q^{-1}b_0^{-1}zy, &\pi(\mathsf{f})(y) &=b_0yz^*,&
\\ \pi(\mathsf{e})(z) &=qb_0^{-1}z^2, &\pi(\mathsf{f})(z) &=-b_0,&
\\ \pi(\mathsf{e})(z^*) &=-q^{-1}b_0^{-1}, &\pi(\mathsf{f})(z^*) &=q^2b_0z^{*2},&
\\ b_0 &\in\mathbb{C},\quad |b_0|^2=q^{-3}.&&&
\end{align*}
In view of Remark \ref{1is}, this set of symmetries is also a part of the series $\mathbf{(1b)}$ distinguished by setting there $a_1=0$, $|a_0|^2=q$.

\item[(iv)] Suppose that $q=\lambda i$ with $\lambda>0$. Then the involution $\mathbf{(D)}$ possesses the compatibility property \eqref{invcomp} with the involution in $\operatorname{Pol}(\mathbb{D})_q$ under a part of the series $\mathbf{(1a)}$ distinguished by setting there $b_1=0$, $|b_0|^2=\lambda^{-3}$. Explicitly, those symmetries are
\begin{align*}
\pi(\mathsf{k})(z) &=q^2z, &\pi(\mathsf{k})(z^*) &=q^{-2}z^*,&
\\ \pi(\mathsf{e})(y) &=q^{-1}b_0^{-1}zy, &\pi(\mathsf{f})(y) &=b_0yz^*,&
\\ \pi(\mathsf{e})(z) &=qb_0^{-1}z^2, &\pi(\mathsf{f})(z) &=-b_0,&
\\ \pi(\mathsf{e})(z^*) &=-q^{-1}b_0^{-1}, &\pi(\mathsf{f})(z^*) &=q^2b_0z^{*2},&
\\ b_0 &\in\mathbb{C},\quad |b_0|^2=\lambda^{-3}.&&&
\end{align*}
In view of Remark \ref{1is}, this set of symmetries is also a part of the series $\mathbf{(1b)}$ distinguished by setting there $a_1=0$, $|a_0|^2=\lambda$.

\item[(v)] Suppose $q=\lambda i$ with $\lambda<0$. Then the involution $\mathbf{(E)}$ possesses the compatibility property \eqref{invcomp} with the involution in $\operatorname{Pol}(\mathbb{D})_q$ under a part of the series $\mathbf{(1a)}$ distinguished by setting there $b_1=0$, $|b_0|^2=-\lambda^{-3}$. Explicitly, those symmetries are
\begin{align*}
\pi(\mathsf{k})(z) &=q^2z, &\pi(\mathsf{k})(z^*) &=q^{-2}z^*,&
\\ \pi(\mathsf{e})(y) &=q^{-1}b_0^{-1}zy, &\pi(\mathsf{f})(y) &=b_0yz^*,&
\\ \pi(\mathsf{e})(z) &=qb_0^{-1}z^2, &\pi(\mathsf{f})(z) &=-b_0,&
\\ \pi(\mathsf{e})(z^*) &=-q^{-1}b_0^{-1}, &\pi(\mathsf{f})(z^*) &=q^2b_0z^{*2},&
\\ b_0 &\in\mathbb{C},\quad |b_0|^2=-\lambda^{-3}.&&&
\end{align*}
In view of Remark \ref{1is}, this set of symmetries is also a part of the series $\mathbf{(1b)}$ distinguished by setting there $a_1=0$, $|a_0|^2=-\lambda$.
\end{description}
The above list is exhaustive. There exist no other $U_q(\mathfrak{sl}_2)$-symmetries on $\operatorname{Pol}(\mathbb{D})_q$ which, under presence of involutions both in $U_q(\mathfrak{sl}_2)$ and in $\operatorname{Pol}(\mathbb{D})_q$, admit the compatibility condition \eqref{invcomp}.
\end{theorem}

{\bf Proof.} The only case which is not covered by Lemma \ref{icex} is {\bf (i)} under the involution $\mathbf{(C)}$. Let us consider this case separately.

Note that, with $\xi\in U_q(\mathfrak{sl}_2)$ being fixed, both l.h.s and r.h.s. of \eqref{invcomp} are anti-linear with respect to $a\in\operatorname{Pol}(\mathbb{D})_q$. Hence it suffices to verify \eqref{invcomp} on a basis. For this purpose, we choose the basis of weight vectors $\left\{\left.z^ky^n,\;y^nz^{*l}\right|\:k,n\ge 0,\,l>0\right\}$ for a symmetry $\pi$. We also consider the basis $\left\{\left.\mathsf{e}^i\mathsf{k}^m\mathsf{f}^j\right|\:i,j\ge 0,\,m\in\mathbb{Z}\right\}$ in $U_q(\mathfrak{sl}_2)$ \cite{KS}, along with an expansion of an arbitrary $\xi\in U_q(\mathfrak{sl}_2)$ with respect to this basis, $\xi=\sum\limits_{i,m,j}c_{i,m,j}\mathsf{e}^i\mathsf{k}^m\mathsf{f}^j$ (the sum is finite). Since both $\pi(\mathsf{e})$ and $\pi(\mathsf{f})$ are identically zero operators, one has in the subcase $\mathbf{(0-)}$
\begin{multline*}
\left(\pi(\xi)\left(z^ky^n\right)\right)^*=
\left(\pi\left(\sum\limits_mc_{0,m,0}\mathsf{k}^m\right)
\left(z^ky^n\right)\right)^*=
\left(\left(\sum\limits_m(-1)^{km}c_{0,m,0}\right)\left(z^ky^n\right)\right)^*
\\ =\left(\sum\limits_m(-1)^{km}\overline{c_{0,m,0}}\right)y^n z^{*k},
\end{multline*}
\begin{multline*}
\left(\pi(\xi)\left(y^nz^{*l}\right)\right)^*=
\left(\pi\left(\sum\limits_mc_{0,m,0}\mathsf{k}^m\right)
\left(y^nz^{*l}\right)\right)^*=
\left(\left(\sum\limits_m(-1)^{lm}c_{0,m,0}\right)
\left(y^nz^{*l}\right)\right)^*
\\ =\left(\sum\limits_m(-1)^{lm}\overline{c_{0,m,0}}\right)z^ly^n,
\end{multline*}
$$
\mathsf{S}(\xi)^*=
\left(\sum\limits_{i,m,j}c_{i,m,j}\mathsf{e}^i\mathsf{k}^m\mathsf{f}^j\right)^*
=\sum\limits_{i,m,j}\overline{c_{i,m,j}}\mathsf{f}^j\mathsf{k}^m\mathsf{e}^i,
$$
$$
\pi(\mathsf{S}(\xi)^*)\left(z^ky^n\right)^*=
\pi\left(\sum\limits_m\overline{c_{0,m,0}}\mathsf{k}^m\right)
\left(y^nz^{*k}\right)=
\left(\sum\limits_m(-1)^{km}\overline{c_{0,m,0}}\right)y^nz^{*k},
$$
$$
\pi(\mathsf{S}(\xi)^*)\left(y^nz^{*l}\right)^*=
\pi\left(\sum\limits_m\overline{c_{0,m,0}}\mathsf{k}^m\right)
\left(z^ly^n\right)=
\left(\sum\limits_m(-1)^{lm}\overline{c_{0,m,0}}\right)z^ly^n,
$$
which establishes \eqref{invcomp}. A similar but even easier argument works also in the subcase $\mathbf{(0+)}$. The claim {\bf (i)} under the involution $\mathbf{(C)}$ is proved.

In all other cases Lemma \ref{icex} is applicable. The latter Lemma allows extraction of suitable subseries satisfying \eqref{invcomp} from the series of symmetries listed explicitly in Sections \ref{trl}, \ref{GJ+}, \ref{GJ-} via verifying \eqref{invcomp} on the generators both in $U_q(\mathfrak{sl}_2)$ and in $\operatorname{Pol}(\mathbb{D})_q$. The verification procedure anticipates calculations which are completely routine and thus left to the reader. \hfill $\blacksquare$

\section*{Acknowledgement}

The author would like to acknowledge the helpful discussions with S. Duplij on the subject of this work.


\begin{thebibliography}{99}
\bibitem{abe} {\it E. Abe}, Hopf Algebras. Cambridge Univ. Press,  Cambridge, 1980.

\bibitem{ale/cha} {\it J. Alev and M. Chamarie}, D\'erivations et
    automorphismes de quelques alg\`ebres quantiques. {\it Comm. Algebra}
    (1992), v. 20, p. 1787 -- 1802.

\bibitem{DHL} {\it S. Duplij, Y. Hong, and F. Li}, $U_q(\mathfrak{sl}_{m+1})$-module algebra structures on the coordinate algebra of a quantum vector space. {\it J. Lie Theory} {\bf 25} (2015), No 2, 327 -- 361.

\bibitem{DS} {\it S. Duplij and S. Sinel'shchikov}, Classification of $U_q(\mathfrak{sl}_2)$-module algebra structures on the quantum plane,   {\it J. Math. Phys., Anal., and Geom.} {\bf 6} (2010), No 4, 406 -- 430.

\bibitem{K} {\it C. Kassel}, Quantum Groups. Springer–Verlag, New York, 1995.

\bibitem{GR} {\it G. Gasper, M. Rahman.} Basic Hypergeometric Series, Cambridge University Press, Cambridge, 1990.

\bibitem{KL} {\it S. Klimek, A. Lesniewski}, A two-parameter quantum deformation of the unit disc, {\it J. Funct. Anal.} {\bf 115} (1993), 1 -- 23.

\bibitem{KS} {\it A. Klimyk, K. Schm{\"u}dgen}. Quantum Groups and Their Representations, Berlin: Springer (1997), 552 pp.

\bibitem{NN} {\it G. Nagy, A. Nica}. On the `quantum disc' and a `non-commutative circle', in: {\it Algebraic Methods on Operator Theory}, R. E. Curto, P. E. T. Jorgensen (eds.), Birkhauser, Boston, 1994, 276 -- 290.

\bibitem{SZ} {\it D. Shklyarov, G. Zhang}, Covariant q-differential operators and unitary highest weight representations for $U_q\mathfrak{su}_{n,n}$, {\it J. Math. Phys.} {\bf 46} (2005), No 6, 24 pp.

\bibitem{S1} {\it S. Sinel'shchikov}, Generic symmetries of the Laurent extension of quantum plane, {\it J. Math. Phys., Anal., and Geom.} {\bf 11} (2015), No 4, 333--358.

\bibitem{S2} {\it S. Sinel'shchikov}, The Laurent extension of quantum plane: a complete list of $U_q(\mathfrak{sl}_2)$-symmetries, {\it SIGMA} {\bf 15} (2019), 038, 33 pp.

\bibitem{SSV1} {\it D. Shklyarov, S. Sinel'shchikov, and L. Vaksman}, q-analogues of some bounded symmetric domains, {\it Czechoslovak Journal of Physics} {\bf 50} (2000), No 1, 175 -- 180.

\bibitem{SSV2} {\it D. Shklyarov, S. Sinel'shchikov, and L. Vaksman}, Geometric realizations for some series of representations of the quantum group $SU_{2,2}$, {\it Math. Phys., Anal., and Geometry} {\bf 8} (2001), No 1, 90 -- 110.

\bibitem{sweedler} {\it M. E. Sweedler}, Hopf Algebras. Benjamin, New York,    1969.

\bibitem{V} {\it L. L. Vaksman}, Quantum bounded symmetric domains, {\it Translations of Mathematical Monographs}, Vol. 238, American Mathematical Society, Providence, RI, 2010.
\end{thebibliography}
\end{document}